\documentclass[11pt]{amsart}
\usepackage{amssymb,amsmath}
\usepackage{graphicx}
\usepackage{amsbsy}
\usepackage{mathrsfs}
\usepackage{color}
\usepackage[active]{srcltx}

\newtheorem{proposition}{Proposition}[section]
\newtheorem{theorem}{Theorem}[section]
\newtheorem{lemma}{Lemma}[section]
\newtheorem{corollary}{Corollary}[section]

\newtheorem{remark}{Remark}[section]
\numberwithin{equation}{section}

\begin{document}

\begin{center}
{\Large \textbf{ Singular Poisson equations on Finsler-Hadamard manifolds}}\\[%
0pt]
\vspace{0.5cm} {\large Csaba Farkas$^*$, Alexandru Krist\'aly$^{**,***,}$%
\footnote{%
Corresponding author. Email address: alexandrukristaly@yahoo.com}, Csaba
Varga$^{****}$}\\[0pt]
\vspace{0.3cm} {\small $^*$ Faculty of Technical and Human Sciences,  Sapientia Hungarian University of Transylvania, 547367 T\^\i rgu Mure\c s, Romania}\\[0pt]
{\small $^{**}$ Department of Economics, Babe\c s-Bolyai University,
400593 Cluj-Napoca, Romania}\\
{\small $^{***}$ Institute of Applied
Mathematics, \'Obuda University, 1034 Budapest, Hungary}\\
{\small $^{****}$ Department of Mathematics, Babe\c s-Bolyai
University, 400084 Cluj-Napoca, Romania}\\[0pt]
\end{center}
\vspace{1cm}

\begin{center} {\it Dedicated to Professor Gheorghe Moro\c sanu on the
occasion of his $65^{th}$ birthday.} \end{center}

\vspace{1cm}

\noindent \textbf{Abstract.} {\footnotesize \noindent In the first
part of the paper we study the reflexivity of Sobolev spaces on
non-compact and not necessarily reversible Finsler manifolds. Then,
by using direct methods in the calculus of variations, we establish
uniqueness, location and rigidity results for singular Poisson
equations involving the Finsler-Laplace operator on Finsler-Hadamard
manifolds having finite reversibility constant.}

\vspace{1cm} \noindent \textbf{Keywords}: 
Poisson equation; Finsler-Hadamard manifold; reversibility constant;
uniformity constant; Finsler-Laplacian operator.\newline


\section{Introduction}

Elliptic problems on Riemannian manifolds have been intensively
studied in the last decades. On one hand, deep achievements have
been done in connection with the famous Yamabe problem on Riemannian
manifolds which can be transformed into an elliptic PDE involving
the  Laplace-Beltrami operator, see Aubin \cite{Aubin} and Hebey
\cite{Hebey}. On the other hand, various anisotropic elliptic
problems are discussed on Minkowski spaces $(\mathbb{R}^n,F)$ where
$F\in C^2(\mathbb{R}^n,[0,\infty))$ is convex and the leading term
is given by the non-linear Finsler-Laplace operator associated with
the Minkowski norm $F$, see Alvino, Ferone, Lions and Trombetti
\cite{AIHP-Lions}, Bellettini and Paolini \cite{bellettini},
Belloni, Ferone and Kawohl \cite{BFK_ZAMP}, \cite{FeKa}, and
references therein. In both classes of problems variational
arguments are applied, the key roles being played by fine properties
of Sobolev spaces  as well as the lower semicontinuity of the energy
functional  associated to the studied problems.

 In order to have a global approach,  the theory of Sobolev spaces has been deeply investigated on metric measure
 spaces, see Ambrosio, Colombo and Di Marino
 \cite{Ambrosioetal}, Cheeger \cite{Cheeger}, and Hajlasz and Koskela
 \cite{Hajlasz}.
In \cite{Ambrosioetal}, the authors proved that if $(X, \mathsf{d})$
is doubling and separable, and $\mathsf{m}$ is finite on bounded
sets, the Sobolev space $W^{1,2}(X,\mathsf{d,m})$ is {\it
reflexive}; here, $W^{1,2}(X,\mathsf{d,m})$ contains functions $u
\in L^2( X, \mathsf{m})$ with finite $2-$relaxed slope endowed by
the norm $u\mapsto \left(\displaystyle\int_X |\nabla u|_{*,2}^2 {\rm
d}\mathsf{m}+\int_X u^2 {\rm d}\mathsf{m}\right)^{1/2}$, where
$|\nabla u|_{*,2}(x)$ denotes the $2-$relaxed slope of $u$ at $x\in
X$. This result clearly applies for differential structures. Indeed,
if $(M,F)$ is a reversible Finsler manifold (in particular, a
Riemannian manifold), then for every $x\in M$ and $u\in
C_0^\infty(M)$,
$$|\nabla u|_{*,2}(x)=\limsup_{z\to x}\frac{|u(z)-u(x)|}{d_F(x,z)}=F^*(x,Du(x)),$$ where $d_{F}$ is the metric function associated with $F$, and $F^*$
is the polar transform of $F,$ see Ohta and Sturm \cite{Ohta-Sturm}.
Consequently, within the class of reversible Finsler manifolds, the
synthetic notion of Sobolev spaces on metric measure spaces (see
\cite{Ambrosioetal} and \cite{Cheeger}) and the analytic notion of
Sobolev spaces on Finsler manifolds  (see Ge and Shen
\cite{Ge-Shen}, and Ohta and Sturm \cite{Ohta-Sturm}) coincide.

Although in the aforementioned works the involved metrics are
symmetric, \textit{asymmetry} is abundant in real life. In order to
describe such phenomena, we put ourselves into the context of not
necessarily reversible Finsler manifolds which model various
Randers-type spaces, including the Matsumoto mountain slope metric,
Finsler-Poincar\'e ball model, etc.; see Bao, Chern and Shen
\cite{BCS}. If $M$ is a connected $n$-dimensional $C^{\infty}$
manifold and $TM=\bigcup_{x \in M}T_{x} M $ is its tangent bundle, the pair $%
(M,F)$ is a \textit{Finsler manifold} if the continuous function
$F:TM\to [0,\infty)$ satisfies the conditions:

(a) $F\in C^{\infty}(TM\setminus\{ 0 \});$

(b) $F(x,ty)=tF(x,y)$ for all $t\geq 0$ and $(x,y)\in TM;$

(c) $g_{ij}(x,y):=\left[ \frac{1}{2}F^{2}(x,y)\right] _{y^{i}y^{j}}$
is positive definite for all $(x,y)\in TM\setminus \{0\}.$\\
If $F(x,ty)=|t|F(x,y)$ for all $t\in\mathbb R$ and $(x,y)\in TM,$ we
say that the Finsler manifold $(M,F)$ is reversible.

Let $(M,F)$ be a Finsler manifold.  Although it is possible to use
an arbitrary measure on $(M,F)$ to define Sobolev spaces (see
\cite{Ohta-Sturm}), here and in the sequel, we shall use the
canonical Hausdorff measure on  $(M,F)$,
$$\mathrm{d}\mathsf{m}=\mathrm{d}V_{F},$$ see Section \ref{sect-2}. Having this measure in our mind, we consider the Sobolev spaces associated with $(M,F)$, see
\cite{Ge-Shen} and  \cite{Ohta-Sturm}. To be more precise, let
\begin{equation*}
W^{1,2}(M,F,\mathsf{m})=\left\{ u\in W_{\mathrm{loc}}^{1,2}(M):\displaystyle%
\int_{M}F^{\ast 2}(x,Du(x))\mathrm{d}\mathsf{m}(x)<+\infty \right\}
,
\end{equation*}%
and $W_{0}^{1,2}(M,F,\mathsf{m})$ be the closure of $C_{0}^{\infty
}(M)$ with respect to the (asymmetric) norm
\begin{equation}\label{Sobolev-norm}
\Vert u\Vert _{F}=\left( \displaystyle\int_{M}F^{\ast 2}(x,Du(x))\mathrm{d}%
\mathsf{m}(x)+\displaystyle\int_{M}u^{2}(x)\mathrm{d}\mathsf{m}(x)\right)
^{1/2}.
\end{equation}%
 Let
\begin{equation*}
r_{F}=\sup_{x\in M}\sup_{\substack{ y\in T_{x}M\setminus \{0\}}}\frac{F(x,y)%
}{F(x,-y)}
\end{equation*}
be the reversibility constant on $(M,F)$. Clearly, $r_F\geq 1$ and
$r_F=1$ if and only if $(M,F)$ is reversible.  Let
\begin{equation*}
F_s(x,y)=\left(\frac{F^2(x,y)+F^2(x,-y)}{2}\right)^{1/2},\ (x,y)\in
TM.
\end{equation*}
It is clear that $(M,F_s)$ is a reversible Finsler manifold, $F_s$
being the {\it symmetrized} Finsler metric associated with $F$. We
notice that the symmetrized Finsler metric associated with $F^*$ may
be different from $F_s^*$, i.e., in general
$2{F_s^*}^2(x,\alpha)\neq {F^*}^2(x,\alpha)+{F^*}^2(x,-\alpha);$
such a concrete case is shown for Randers metrics, see
(\ref{szimmetrikus-randers}).

Our first result reads as follows:

\begin{theorem}\label{thm-1} Let $(M,F)$ be a complete, $n-$dimensional Finsler
manifold such that $r_F<+\infty$.  Then
$(W_{0}^{1,2}(M,F,\mathsf{m}),\|\cdot\|_{F_s})$ is a reflexive
Banach space, while the norm $\|\cdot\|_{F_s}$ and the asymmetric
norm $\|\cdot\|_{F}$ are equivalent. In particular,
\begin{equation}\label{norm-equivalent}
\left(\frac{1+r_F^{2}}{2}\right)^{-{1}/{2}}\|u\|_{F}\leq
\|u\|_{F_s}\leq
\left(\frac{1+r_F^{-2}}{2}\right)^{-{1}/{2}}\|u\|_{F},\ \forall u\in
W_{0}^{1,2}(M,F,\mathsf{m}).
\end{equation}
\end{theorem}
For sake of clarity, we notice that the norm $\|\cdot\|_{F_s}$ is
considered also with respect to the Hausdorff measure
$\mathrm{d}\mathsf{m}=\mathrm{d}V_{F}$ (and not with
$\mathrm{d}V_{F_s}$), i.e.,
\begin{equation}\label{Sobolev-norm-Fs}
\Vert u\Vert _{F_s}=\left( \displaystyle\int_{M}F_s^{\ast 2}(x,Du(x))\mathrm{d}%
\mathsf{m}(x)+\displaystyle\int_{M}u^{2}(x)\mathrm{d}\mathsf{m}(x)\right)
^{1/2}.
\end{equation}
 Some remarks are in order concerning Theorem \ref{thm-1}.
\begin{remark}\rm
(i) We emphasize that Theorem \ref{thm-1} is sharp. Indeed, let us
consider the two-dimensional Finsler-Poincar\'e model $(B^2(0,2),F)$
which is a forward (but not backward) complete Finsler manifold of
Randers-type having the reversibility constant $r_F=+\infty$, see
Section \ref{sect-3}. In this framework, we shall construct a
function $u\in W_0^{1,2}(B^2(0,2),F,\mathsf{m})$ such that $
 -u\notin W_0^{1,2}(B^2(0,2),F,\mathsf{m});$ in other words,
 $W_0^{1,2}(B^2(0,2),F,\mathsf{m})$  {\it does not have} a vector space structure, and the norm $\|\cdot\|_{F_s}$ and the asymmetric norm
$\|\cdot\|_{F}$ are not equivalent.  A similar pathological
situation has been already pointed out by Krist\'{a}ly and Rudas
\cite{Kristaly-NA} for a Funk-type metric on
the open unit ball of $\mathbb{R}^{n}$. 

(ii) It is clear that $r_F<+\infty$ whenever $(M,F)$ is a
\textit{compact} Finsler manifold. Thus, the reflexivity of
$W_{0}^{1,2}(M,F,\mathsf{m})$ in \cite{Ge-Shen} and
\cite{Ohta-Sturm} immediately follows from Theorem \ref{thm-1}.

(iii) We believe that $(W_{0}^{1,2}(M,F,\mathsf{m}),\|\cdot\|_{F})$
is a reflexive, complete asymmetric vector space; a possible proof
requires a long series of arguments from functional analysis for
asymmetric normed spaces, see Cobza\c s \cite{Cobzas}. However, the
statement of Theorem \ref{thm-1} is enough for our purposes.
\end{remark}

In the second part we consider that $(M,F)$ is an $n-$dimensional
Finsler-Hadamard manifold (i.e., simply connected, complete with
non-positive flag curvature), $n\geq 3$, having its uniformity constant $%
l_{F}>0$ (which implies in particular that $r_{F}<+\infty $), see
Section \ref{sect-2}. We shall study the model singular Poisson
equation
\begin{equation*}
\ \left\{
\begin{array}{lll}
\boldsymbol{\Delta }(-u)-\mu \frac{u}{d_{F}^{2}(x_{0},x)}=1 & \mbox{in} &
\Omega ; \\
u=0 & \mbox{on} & \partial \Omega,
\end{array}%
\right. \eqno{({\mathcal P}^{\mu}_\Omega)}
\end{equation*}%
where $\boldsymbol{\Delta }$ denotes the Finsler-Laplace operator on $(M,F)$%
,  $x_{0}\in \Omega $ is fixed, $\mu \geq 0$ is a parameter, and
$\Omega \subset M$ is an open and bounded domain with sufficiently
smooth boundary. We prove that the singular energy functional
associated with problem $({\mathcal{P}}_{\Omega }^{\mu })$ is {\it
strictly
convex} on $W_{0}^{1,2}(\Omega,F,\mathsf{m} )$ whenever $\mu \in \lbrack 0,l_{F}r_{F}^{-2}%
\overline{\mu })$, see Theorem \ref{energy-convex}; here, $\overline{\mu }=%
\frac{(n-2)^{2}}{4}$ is the optimal Hardy constant. By exploiting
Theorem \ref{thm-1}, a comparison principle for the Finsler-Laplace
operator and well known arguments from calculus of variations, we
prove (see also  Theorem \ref{unique-theo}):

\begin{theorem}\label{uniq-elso} Problem $({\mathcal{%
P}}_{\Omega }^{\mu })$ has a unique, non-negative weak solution whenever $%
\mu \in \lbrack 0,l_{F}r_{F}^{-2}\overline{\mu })$.
\end{theorem}
Having the uniqueness theorem in our mind, we focus our attention to
geometric rigidities related to the Poisson equation $({\mathcal
P}^{\mu}_\Omega)$. To do this, let $c\leq 0$ and the function ${\bf
ct}_{c}:(0,\infty)\to \mathbb R$ defined by
$${\bf ct}_{c}(r)=\left\{
  \begin{array}{lll}
    \frac{1}{r}
    & \hbox{if} &  {c}=0, \\

  \sqrt{-c}\coth(\sqrt{-c}r) & \hbox{if} & {c}<0.
  \end{array}\right.$$
For every  $\mu \in \lbrack 0,%
\overline{\mu })$, $\rho>0$ and $c\leq 0$, we consider the
 ordinary differential equation
\begin{equation*}
\left\{
\begin{array}{ll}
f^{\prime \prime }(r)+(n-1) f^{\prime }(r){\bf ct}_{c}(r)+\mu
\frac{f(r)}{r^{2}}+1=0, \ r\in (0,\rho ], \\
f(\rho )=0,\ \displaystyle\int_0^\rho f'(r)^2r^{n-1}{\rm d}r<\infty.
\end{array}%
\right. \eqno{({\mathcal Q}^{\mu}_{c,\rho})}
\end{equation*}
We shall show that $({\mathcal Q}^{\mu}_{c,\rho})$ has a unique,
non-negative non-increasing solution  $\sigma _{\mu ,\rho ,c}\in
C^\infty(0,\rho)$, see Proposition \ref{prop-ODE}. Although we are
not able to solve explicitly $({\mathcal Q}^{\mu}_{c,\rho})$, in
some particular cases we have its solution; namely,
$$\sigma _{\mu ,\rho ,c}(r)=\left\{
  \begin{array}{lll}
   \frac{1}{\mu +2n}\left( \rho
^{2}\left( \frac{r}{\rho }\right) ^{-\sqrt{\overline{\mu }}+\sqrt{\overline{%
\mu }-\mu }}-{r}^{2}\right)
    & \hbox{if} &  {c}=0;
    \\
  \displaystyle\int_{r}^{\rho}\sinh (\sqrt{-c}s)^{-n+1}\int_{0}^{s}\sinh (\sqrt{-c}%
t)^{n-1}{\rm d}t{\rm d}s & \hbox{if} & {c}<0\ {\rm and}\ \mu=0;\\
H(\sqrt{\overline{%
\mu }-\mu },\rho)\frac{\sqrt{r}\sinh(\rho)I_{\sqrt{\overline{%
\mu }-\mu }}(r)}{\sqrt{\rho}\sinh(r)I_{\sqrt{\overline{%
\mu }-\mu }}(\rho)}-H(\sqrt{\overline{%
\mu }-\mu },r) & \hbox{if} & {c}=-1,\ n=3\ {\rm and}\ \mu\in
[0,\frac{1}{4}),
  \end{array}\right.$$
  where $H:(0,\frac{1}{2}]\times (0,\infty)\to \mathbb R$ is given by
\begin{eqnarray*}
  H(\nu,r) &=& \frac{2\nu}{({25}-4\nu^2)\sin(\nu\pi)\Gamma(\nu)\sinh(r)}\times \\
   &&\times \left\{\left(5-2\nu\right) {_3F{_4}}\left(\left[\frac{3}{4}+\frac{\nu}{2},\frac{5}{4}+\frac{\nu}{2},\frac{5}{4}+\frac{\nu}{2}\right];
   \left[\frac{3}{2},1+{\nu},\frac{3}{2}+{\nu},\frac{9}{4}+\frac{\nu}{2}\right],r^2\right) \right.\times\\
   &&\ \ \ \ \ \ \times\left(2^{\nu-2}\sin(\nu\pi)K_\nu(r)+2^{-\nu-1}\pi I_\nu(r)\right)r^{3+\nu}- \\
   &&\ \ \ \ \ \ \left. -\nu(5+2\nu)2^{\nu-1}{_3F{_4}}\left(\left[\frac{3}{4}-\frac{\nu}{2},\frac{5}{4}-\frac{\nu}{2},\frac{5}{4}-\frac{\nu}{2}\right];
   \left[\frac{3}{2},1-{\nu},\frac{3}{2}-{\nu},\frac{9}{4}-\frac{\nu}{2}\right],r^2\right)\right.\times\\
   &&\ \ \ \ \ \ \ \ \ \bigg. \times\Gamma(\nu)^2\sin(\nu\pi)I_\nu(r)r^{3-\nu}  \displaystyle
   \bigg\};
\end{eqnarray*}
here, $I_\nu$ and $K_\nu$ are the  modified Bessel functions of the
first and second kinds of order $\nu,$ while $_3F_4$ denotes the
generalized hypergeometric function.

Before to state our rigidity results, we need
  some notations:
${\bf K}\leq c$ (resp. $c\leq {\bf K}$, resp. ${\bf K}= c$) means
that the flag curvature ${\bf K}(\mathcal{S};v)$ on $(M,F)$ is
bounded from above by $c\in \mathbb R$ (resp. bounded from below by
$c$, resp. equal to $c$) for any choice of parameters $\mathcal S$
and $v$; ${\bf S}=0$ means that $(M,F)$ has vanishing mean
covariation; $B^{+}(x_{0},\rho )$ denotes the open forward metric
 ball with center $x_0\in M$ and radius $\rho>0;$ for details, see
Section \ref{sect-2}.

\begin{theorem} {\rm (Local estimate via curvature)}
\label{estimate_with_constant_flag} Let $(M,F)$ be an $n-$dimensional $%
(n\geq 3)$ Finsler-Hadamard manifold with $\mathbf{S}=0$ and
$l_{F}>0,$ and $\Omega\subset M$ be an open bounded domain.  Let
$\mu\in [0,l_{F} r_{F}^{-2} \overline \mu)$ and  $x_{0}\in \Omega$
be fixed. If $c_{1}\leq
\mathbf{K}\leq c_{2}\leq 0,$ then the unique weak solution $u$ of problem $({%
\mathcal{P}}_{\Omega}^{\mu })$ verifies the  inequalities%
\begin{equation*}
\sigma _{\mu ,\rho_1 ,c_{1}}(d_{F}(x_{0},x))\leq u(x)\leq \sigma
_{\mu ,\rho_2 ,c_{2}}(d_{F}(x_{0},x))\text{ for a.e.
 }x\in
B^{+}(x_{0},\rho_1 ),
\end{equation*}
where $\rho_1=\sup\{\rho>0:B^+(x_0,\rho)\subset \Omega\}$ and
$\rho_2=\inf\{\rho>0:\Omega\subset B^+(x_0,\rho) \}.$

In particular, if $\mathbf{K}=c\leq 0$ and $\Omega=B^+(x_0,\rho)$
for some $\rho>0,$ then $\sigma _{\mu ,\rho
,c}(d_{F}(x_{0},\cdot))$ is the unique weak solution of problem $({%
\mathcal{P}}_{B^{+}(x_{0},\rho )}^{\mu })$, being also a pointwise
solution in $B^+(x_0,\rho)\setminus \{x_0\}.$
\end{theorem}

\noindent A kind of converse statement of Theorem
\ref{estimate_with_constant_flag} can read as follows.

\begin{theorem}{\rm (Radial curvature rigidity)}
\label{rigidity} Let $(M,F)$ be an $n-$dimensional $(n\geq 3)$
Finsler-Hadamard manifold with $\mathbf{S}=0$, $l_{F}>0$ and
$\mathbf{K}\leq c\leq 0.$ Let $\mu\in [0,l_{F} r_{F}^{-2} \overline
\mu)$ and  $x_{0}\in M$ be fixed. If the function $\sigma _{\mu
,\rho ,c}(d_{F}(x_{0},\cdot))$ is the unique pointwise solution of
$({\mathcal{P}}_{B^{+}(x_{0},\rho )}^{\mu })$ in
$B^+(x_0,\rho)\setminus \{x_0\}$ for some $\rho
>0,$ then  ${\bf K}(\cdot;\dot \gamma_{x_0,y}(t))=c$ for every $t\in
[0,\rho)$ and $y\in T_{x_0}M\setminus \{0\}$, where $\gamma_{x_0,y}$
is the constant speed geodesic with $\gamma_{x_0,y}(0)=x_0$ and
$\dot \gamma_{x_0,y}(0)=y.$
\end{theorem}

In the generic Finsler setting,  the conclusion of Theorem
\ref{rigidity} does not imply necessarily that the flag curvature
${\bf K}$ is constant. Indeed, we just stated that the flag
curvature is {\it radially} constant  with respect to $x_0\in M$,
i.e., along geodesics emanating from the point $x_0$ where the
flag-poles are the velocities of the geodesics. However, when
$(M,F)=(M,g)$ is a {\it Riemannian} manifold of Hadamard type (thus
the flag curvature and sectional curvature coincide and the notion
of the flag loses its meaning), Theorems
\ref{estimate_with_constant_flag}\&\ref{rigidity} and the
classification of Riemannian space forms (see do Carmo \cite[Theorem
4.1]{doCarmo}) provide a {\it characterization} of the Euclidean and
hyperbolic spaces up to isometries via the shape of solutions to the
Poisson equation $({\mathcal P}^{\mu}_\Omega)$:

\begin{corollary}\label{corollary} {\rm (Space forms vs. Poisson equation)} Let
$(M,g)$ be a Riemannian-Hadamard manifold with sectional curvature
bounded above by $c\leq 0.$ Then the following statements are
equivalent:

\begin{itemize}
\item[(a)] For some $\mu \in [0,\overline{\mu })$ and $x_{0}\in M$, the
function $\sigma _{\mu ,\rho,c }(d_{F}(x_{0},\cdot))$ is the unique
pointwise solution of the Poisson equation
$({\mathcal{P}}_{B(x_{0},\rho )}^{\mu })$ in $B(x_0,\rho)\setminus
\{x_0\}$ for every $\rho >0;$

\item[(b)] $(M,g)$ is isometric to the $n-$dimensional space form
with  curvature $c.$
\end{itemize}
\end{corollary}

 A full classification of
Finslerian space forms (i.e., the flag curvature is constant) is not
available; however, the following characterization can be provided
on {\it Berwald spaces}:

\begin{theorem}{\rm (Full curvature rigidity)}
\label{rigid_Minkowski} Let $(M,F)$ be an $n-$dimensional $(n\geq
3)$ Finsler-Hadamard manifold of Berwald type with $l_{F}>0.$ Then
the following statements are equivalent:

\begin{itemize}
\item[(a)] For every $\mu \in \lbrack 0,l_{F}r_{F}^{-2}\overline{\mu })$ and $x_{0}\in M$, the
function $\sigma _{\mu ,\rho,0 }(d_{F}(x_{0},\cdot))$ is the unique
pointwise solution of the Poisson equation
$({\mathcal{P}}_{B^{+}(x_{0},\rho )}^{\mu })$ in
$B^+(x_0,\rho)\setminus \{x_0\}$ for every $\rho >0;$

\item[(b)] For some $\mu \in \lbrack 0,l_{F}r_{F}^{-2}\overline{\mu })$ and $x_{0}\in M$, the
function $\sigma _{\mu ,\rho,0 }(d_{F}(x_{0},\cdot))$ is the unique
pointwise solution of the Poisson equation
$({\mathcal{P}}_{B^{+}(x_{0},\rho )}^{\mu })$ in
$B^+(x_0,\rho)\setminus \{x_0\}$ for every $\rho >0;$

\item[(c)] $(M,F)$ is isometric to an $n-$dimensional Minkowski space.
\end{itemize}
\end{theorem}

\section{Preliminaries: elements from Finsler
geometry}\label{sect-2}

\textit{2.1. Finsler manifolds, geodesics, flag curvature, mean
covariation, volume element.} Let $(M,F)$ be a Finsler manifold
(i.e., (a)-(c) hold from the Introduction). If
$g_{ij}(x)=g_{ij}(x,y)$ is independent of $y$ then $(M,F)=(M,g)$ is
called a \textit{Riemannian manifold}. A \textit{Minkowski space}
consists
of a finite dimensional vector space $V$ (usually, identified with $\mathbb{R%
}^{n}$) and a Minkowski norm which induces a Finsler metric on $V$ by
translation, i.e., $F(x,y)$ is independent on the base point $x$; in such
cases we often write $F(y)$ instead of $F(x,y)$.

A specific non-reversible Finsler structure is provided by \textit{Randers
metrics} which will serve to us as a model case. To be more precise, on a
manifold $M$ we introduce the Finsler structure $F:TM\rightarrow \lbrack
0,\infty )$ defined by
\begin{equation}
F(x,y)=\sqrt{h_{x}(y,y)}+\beta _{x}(y),\ (x,y)\in TM,
\label{Randers-metrika}
\end{equation}%
where $h$ is a Riemannian metric on $M$, $\beta $ is an 1-form on $M$, and
we assume that
\begin{equation*}
\Vert \beta \Vert _{h}(x)=\sqrt{h_{x}^{\ast }(\beta _{x},\beta _{x})}<1,\
\forall x\in M.
\end{equation*}%
Here, the co-metric $h_{x}^{\ast }$ can be identified by $h_{x}^{-1}$, the
inverse of the symmetric, positive definite matrix $h_{x}$. Clearly, the
Randers space  $(M,F)$ in (\ref{Randers-metrika}) is symmetric if and only if $%
\beta =0$. Note that Randers metrics appear in the study of the
electromagnetic field of the physical space-time in general relativity, see
Randers \cite{Randers}. Moreover, a deep result of Bao, Robles and Shen \cite%
{BRS} shows that a Finsler metric is of Randers type if and only if it is a
solution of the Zermelo navigation problem on a Riemannian manifold.

Let $\pi ^{*}TM$ be the pull-back bundle of the tangent bundle $TM$
generated by the natural projection $\pi:TM\setminus\{ 0 \}\to M,$ see Bao,
Chern and Shen \cite[p. 28]{BCS}. The vectors of the pull-back bundle $\pi
^{*}TM$ are denoted by $(v;w)$ with $(x,y)=v\in TM\setminus\{ 0 \}$ and $%
w\in T_xM.$ For simplicity, let $\partial_i|_v=(v;\partial/\partial x^i|_x)$
be the natural local basis for $\pi ^{*}TM$, where $v\in T_xM.$ One can
introduce on $\pi ^{*}TM$ the \textit{fundamental tensor} $g$
by
\begin{equation}  \label{funda-Cartan-tensors}
g_v:=g(\partial_i|_v,\partial_j|_v)=g_{ij}(x,y)
\end{equation}
where $v=y^i{(\partial}/{\partial x^i})|_x.$ Unlike the Levi-Civita
connection
\index{Levi-Civita connection} in the Riemannian case, there is no unique
natural connection in the Finsler geometry. Among these connections on the
pull-back bundle $\pi ^{*}TM,$ we choose a torsion free and almost
metric-compatible linear connection on $\pi ^{*}TM$, the so-called \textit{%
Chern connection}, see Bao, Chern and Shen \cite[Theorem 2.4.1]{BCS}. The
coefficients of the Chern connection are denoted by $\Gamma_{jk}^{i}$, which
are instead of the well known Christoffel symbols from Riemannian geometry.
A Finsler manifold is of \textit{Berwald type} if the coefficients $%
\Gamma_{ij}^{k}(x,y)$ in natural coordinates are independent of $y$. It is
clear that Riemannian manifolds and $($locally$)$ Minkowski spaces are
Berwald spaces. The Chern connection induces on $\pi ^{*}TM$ the \textit{%
curvature tensor} $R$, see Bao, Chern and Shen \cite[Chapter 3]{BCS}. By
means of the connection, we also have the \textit{covariant derivative} $%
D_vu $ of a vector field $u$ in the direction $v\in T_xM$ with reference
vector $v $. A vector field $u=u(t)$ along a curve $\sigma$ is \textit{%
parallel} if $D_{\dot \sigma}u=0.$ A $C^\infty$ curve
$\sigma:[0,a]\to M$ is a \textit{geodesic} if $D_{\dot \sigma}{\dot
\sigma}=0.$ Geodesics are considered to be parametrized
proportionally to arc-length. The Finsler manifold is
\textit{forward} \textrm{(resp.} {\it backward$)$ complete} if every
geodesic segment $\sigma:[0,a]\to M$ can be extended to $[0,\infty)$
(resp. to $(-\infty,a]$). $(M,F)$ is \textit{complete} if it is both
forward and backward complete.

Let $u,v\in T_xM$ be two non-collinear vectors and $\mathcal{S}=\mathrm{span}%
\{u,v\}\subset T_xM$. By means of the curvature tensor $R$, the \textit{flag
curvature} of the flag $\{\mathcal{S},v\}$ is defined by
\begin{equation}  \label{ref-flag}
\mathbf{K}(\mathcal{S};v) =%
\frac{g_v(R(U,V)V, U)}{g_v(V,V) g_v(U,U) - g_v(U,V)^{2}},
\end{equation}
where $U=(v;u),V=(v;v)\in \pi^*TM.$ If $(M,F)$ is Riemannian, the flag
curvature reduces to the well known sectional curvature. If $\mathbf{K}(%
\mathcal{S};v)\leq 0$ for every choice of $U$ and $V$, we say that $(M,F)$
has \textit{non-positive flag curvature}, and we denote by $\mathbf{K}\leq 0$%
. $(M,F)$ is a \textit{Finsler-Hadamard manifold} if it is simply connected,
forward complete with $\mathbf{K}\leq 0$.

Let $\sigma: [0,r]\to M$ be a piecewise $C^{\infty}$ curve. The value $%
L_F(\sigma)= \displaystyle\int_{0}^{r} F(\sigma(t), \dot\sigma(t))\,{\text d}%
t $ denotes the \textit{integral length} of $\sigma.$ For $x_1,x_2\in M$,
denote by $\Lambda(x_1,x_2)$ the set of all piecewise $C^{\infty}$ curves $%
\sigma:[0,r]\to M$ such that $\sigma(0)=x_1$ and $\sigma(r)=x_2$. Define the
\textit{distance function} $d_{F}: M\times M \to[0,\infty)$ by
\begin{equation}  \label{quasi-metric}
d_{F}(x_1,x_2) = \inf_{\sigma\in\Lambda(x_1,x_2)} L_F(\sigma).
\end{equation}
One clearly has that $d_{F}(x_1,x_2) =0$ if and only if $x_1=x_2,$ and $d_F$
verifies the triangle inequality. The open \textit{forward} (resp. \textit{%
backward}) \textit{metric ball} with center $x_0\in M$ and radius $\rho>0$
is defined by $B^+(x_0,\rho)=\{x\in M:d_F(x_0,x)< \rho \}$ (resp. $%
B^-(x_0,\rho)=\{x\in M:d_F(x,x_0)< \rho \}$). In particular, when $(M,F)=(%
\mathbb{R}^n,F)$ is a Minkowski space, one has $d_F(x_1,x_2)=F(x_2-x_1).$

Let $\{{\partial}/{\partial x^i} \}_{i=1,...,n}$ be a local basis for the
tangent bundle $TM,$ and $\{\mathrm{d}x^i \}_{i=1,...,n}$ be its dual basis
for $T^*M.$ Let $B_x(1)=\{y=(y^i):F(x,y^i \partial/\partial x^i)< 1\}\subset \mathbb R^n$. The \textit{Hausdorff volume form} ${\text d}{\mathsf m}={\text d%
}V_F$ on $(M,F)$ is defined by
\begin{equation}  \label{volume-form}
{\text d}{\mathsf m}(x)={\text d}V_F(x)=\sigma_F(x){\text
d}x^1\wedge...\wedge {\text d}x^n,
\end{equation}
where $\sigma_F(x)=\frac{\omega_n}{\mathrm{Vol}(B_x(1))}$. Hereafter, Vol$%
(S) $ and $\omega_n$ denote the Euclidean volume of the set $S\subset
\mathbb{R}^n$ and the $n-$dimensional unit ball, respectively. The \textit{%
Finslerian-volume} of an open set $S\subset M$ is
Vol$_F(S)=\displaystyle\int_S {\text d}{\mathsf m}(x)$.

Let $\{e_i\}_{i=1,...,n}$ be a basis for $T_xM$ and $g_{{ij}}^v=g_v(e_i,e_j)$%
. The \textit{mean distortion} $\mu:TM\setminus \{0\}\to (0,\infty)$
is defined by
$\mu(v)=\frac{\sqrt{\mathrm{det}(g_{ij}^v)}}{\sigma_F}$. The
\textit{mean covariation} $\mathbf{S}:TM\setminus\{0\}\to
\mathbb{R}$ is defined by
\begin{equation*}
\mathbf{S}(x,v)=\frac{\mathrm{d}}{\mathrm{d}t}(\ln \mu(\dot\sigma_v(t)))\big|%
_{t=0},
\end{equation*}
where $\sigma_v$ is the geodesic such that $\sigma_v(0)=x$ and $\dot
\sigma_v(0)=v.$ 
We say that $(M,F)$ has \textit{vanishing mean covariation} if $\mathbf{S}%
(x,v)= 0$ for every $(x,v)\in TM$, and we denote by $\mathbf{S}=0$.
We
notice that any Berwald space has vanishing mean covariation, see Shen \cite%
{Shen-volume}.

For any $c\leq 0,$ we introduce
$$V_{c,n}(\rho)=n\omega_n\int_0^\rho {\bf s}_c(t)^{n-1}{\rm d}t,$$
where ${\bf s}_c$ denotes the unique solution of $y''+cy=0$ with
$y(0)=0$ and $y'(0)=1,$ i.e.,
$${\bf s}_{c}(r)=\left\{
  \begin{array}{lll}
    r
    & \hbox{if} &  {c}=0, \\
  \frac{\sinh(\sqrt{-c}r)}{\sqrt{-c}} & \hbox{if} & {c}<0.
  \end{array}\right.$$
In general, one has for every $x\in M$ that
\begin{equation}  \label{volume-comp-nullaban}
\lim_{\rho\to 0^+}\frac{\mathrm{Vol}_F(B^+(x,\rho))}{V_{c,n}(\rho)}%
=\lim_{\rho\to
0^+}\frac{\mathrm{Vol}_F(B^-(x,\rho))}{V_{c,n}(\rho)}=1.
\end{equation}

When $(\mathbb{R}^n,F)$ is a Minkowski space, then on account of (\ref%
{volume-form}), Vol$_F(B^+(x,\rho))=\omega_n\rho ^n$ for every $\rho>0$ and $%
x\in \mathbb{R}^n$, and $\sigma_F(x)=$constant. If $F$ is the
Randers metric of the form (\ref{Randers-metrika}) on a manifold
$M$, then
\begin{equation}  \label{randers-volume}
{\text d}V_F(x)=\left(1-\|\beta\|^2_h(x)\right)^\frac{n+1}{2}{\text d}V_h(x),
\end{equation}
where ${\text d}V_h(x)$ denotes the canonical Riemannian volume form of $h$ on $%
M.$

We shall use a Bishop-Gromov volume comparison result; on account of
Shen \cite{Shen-volume}, Wu and Xin \cite[Theorems 6.1 \&
6.3]{Wu-Xin} and Zhao and Shen \cite[Theorem 3.6]{Zhao-Shen}, we
recall the following version:

\begin{theorem}
\label{comparison-volume}{\rm [Volume comparison]} Let $(M,F)$ be an $n-$%
di\-men\-sional Finsler-Hadamard manifold with $\mathbf{S}=0$, ${\bf
K}\leq c\leq 0$ and $x\in M$ fixed. Then the function $$\rho\mapsto
\frac{\mathrm{Vol}_F(B^+(x,\rho))}{V_{c,n}(\rho)}, \ \rho>0,$$ is
non-decreasing.
In particular, from {\rm \textrm{(\ref{volume-comp-nullaban})}} we
have
\begin{equation}  \label{volume-comp-altalanos-0}
{\mathrm{Vol}_F(B^+(x,\rho))}\geq V_{c,n}(\rho)\ {for\ all}\ \rho>0.
\end{equation}
If equality holds in {\rm \textrm{(\ref{volume-comp-altalanos-0})}}
for some $\rho_0>0$, then ${\bf K}(\cdot;\dot \gamma_y(t))=c$ for
every $t\in [0,\rho_0)$ and $y\in T_xM$ with $F(x,y)=1$, where
$\gamma_y$ is the constant speed geodesic with $\gamma_y(0)=x$ and
$\dot \gamma_y(0)=y.$
\end{theorem}


\textit{2.2. Polar and Legendre transforms.} We consider the \textit{polar
transform} (or, co-metric) of $F$, defined for every $(x,\alpha)\in T^*M$ by
\begin{equation}  \label{polar-transform}
F^*(x,\alpha)=\sup_{y\in T_xM\setminus \{0\}}\frac{\alpha(y)}{F(x,y)}.
\end{equation}
Note that for every $x\in M$, the function $F^*(x,\cdot)$ is a Minkowski
norm on $T_x^*M.$ Since $F^{*2}(x,\cdot)$ is twice differentiable on $%
T_x^*M\setminus \{0\}$, we consider the matrix $g_{ij}^*(x,\alpha):=[%
\frac12F^{*2}(x,\alpha)]_{\alpha^{i}\alpha^{j}}$ for every $\alpha=%
\displaystyle\sum_{i=1}^n\alpha^i\mathrm{d}x^i\in T_x^*M\setminus \{0\}$ in
a local coordinate system $(x^i).$

In particular, if $(M,F)$ is a Randers space of the form
(\ref{Randers-metrika}), then
\begin{equation}  \label{Randers-dual}
F^*(x,\alpha)=\frac{\sqrt{h_x^{*2}(\alpha,\beta)+(1-\|\beta\|_h^2(x))\|%
\alpha\|_h^2(x)}-h_x^{*}(\alpha,\beta)}{1-\|\beta\|_h^2(x)},\ (x,\alpha)\in
T^*M,
\end{equation}
where $h_x^*$ denotes the co-metric acting on $T_x^*M$ associated to
the Riemannian metric $h.$ Moreover, the symmetrized Finsler metric
and its polar transform associated with the Randers metric
(\ref{Randers-metrika}) is
\begin{equation}\label{szimmetrikus-randers}
    F_s(x,y)=\sqrt{h_x(y,y)+\beta _{x}^2(y)},\
    F_s^*(x,\alpha)=\sqrt{\|\alpha\|_h^2(x)-\frac{h_x^{*2}(\alpha,\beta)}{1+\|\beta\|_h^2(x)}}.
\end{equation}

The \textit{Legendre transform} $J^*:T^*M\to TM$ associates to each element $%
\alpha\in T_x^*M$ the unique maximizer on $T_xM$ of the map $y\mapsto
\alpha(y)-\frac{1}{2}F^2(x,y)$. This element can also be interpreted as the
unique vector $y\in T_xM$ with the properties
\begin{equation}  \label{ohta-sturm-dual}
F(x,y)=F^*(x,\alpha)\ \mathrm{and}\ \alpha(y)=F(x,y)F^*(x,\alpha).
\end{equation}

In particular, if $\alpha=\sum_{i=1}^n\alpha^i \mathrm{d}x^i\in
T_x^*M$, one has
\begin{equation}  \label{j-csillag}
J^*(x,\alpha)=\sum_{i=1}^n\frac{\partial}{\partial \alpha_i}\left(\frac{1}{2}%
F^{*2}(x,\alpha)\right)\frac{\partial}{\partial x^i}.
\end{equation}

\textit{2.3. Derivatives, Finsler-Laplace operator.} Let $u:M\to \mathbb{R}$
be a differentiable function in the distributional sense. The \textit{%
gradient} of $u$ is defined by
\begin{equation}  \label{grad-deriv}
\boldsymbol{\nabla} u(x)=J^*(x,Du(x)),
\end{equation}
where $Du(x)\in T_x^*M$ denotes the (distributional) \textit{derivative} of $%
u$ at $x\in M.$ In local coordinates, one has
\begin{equation}  \label{derivalt-local}
Du(x)=\sum_{i=1}^n \frac{\partial u}{\partial x^i}(x)\mathrm{d}x^i,
\end{equation}
\begin{equation*}
\boldsymbol{\nabla} u(x)=\sum_{i,j=1}^n g_{ij}^*(x,Du(x))\frac{\partial u}{%
\partial x^i}(x)\frac{\partial}{\partial x^j}.
\end{equation*}

In general, $u\mapsto\boldsymbol{\nabla} u $ is not linear. If $x_0\in M$ is
fixed, then due to Ohta and Sturm \cite{Ohta-Sturm}, one has {\
\begin{equation}  \label{tavolsag-derivalt}
F^*(x,D d_F(x_0,x))=F(x,\boldsymbol{\nabla} d_F(x_0,x))=D d_F(x_0,x)(%
\boldsymbol{\nabla} d_F(x_0,x))=1\ \mathrm{for\ a.e.}\ x\in M.
\end{equation}%
} Let $X$ be a vector field on $M$. In a local coordinate system $(x^i)$, on
account of (\ref{volume-form}), the \textit{divergence} is defined by div$%
(X)=\frac{1}{\sigma_F}\frac{\partial}{\partial x^i}(\sigma_F X^i).$ The
\textit{Finsler-Laplace operator}
\begin{equation*}
\boldsymbol{\Delta} u=\mathrm{div}(\boldsymbol{\nabla} u)
\end{equation*}
acts on $W^{1,2}_{\mathrm{loc}}(M)$ and for every $v\in C_0^\infty(M)$,
\begin{equation}  \label{Green}
\int_M v\boldsymbol{\Delta} u {\text d}{\mathsf m}(x)=-\int_M
Dv(\boldsymbol{\nabla} u){\text d}{\mathsf m}(x),
\end{equation}
see  Ohta and Sturm \cite{Ohta-Sturm} and Shen
\cite{Shen-monograph}. Note that in general $\boldsymbol{\Delta}
(-u)\neq -\boldsymbol{\Delta} u,$
unless $(M,F)$ is  reversible. In particular, for a Riemannian manifold $%
(M,F)=(M,g)$ the Finsler-Laplace operator is the usual Laplace-Beltrami operator $%
\boldsymbol{\Delta}u=\Delta_gu,$ while for a Minkowski space $(\mathbb{R}%
^n,F)$, by using (\ref{ohta-sturm-dual}), $\boldsymbol{\Delta}u=\Delta_Fu=%
\mathrm{div} (F^*(Du)\nabla F^*(Du))=\mathrm{div} (F(\nabla u)\nabla
F(\nabla u))$ is precisely the Finsler-Laplace operator considered
by Cianchi and Salani \cite{Cianchi-Salani}, Ferone and Kawohl
\cite{FeKa}, Wang and Xia \cite{WX, WX-2}, and references therein.
We shall use the following result from Wu and Xin \cite{Wu-Xin}:

\begin{theorem}
\label{comparison-lagrange}{\rm [Laplacian comparison]} Let $(M,F)$ be an $n-$%
di\-men\-sional Finsler-Hadamard manifold with $\mathbf{S}=0$. Let
$x_0\in M$ and $c\leq 0.$ Then the  following statements hold:
\begin{itemize}
  \item[{\rm (a)}] If ${\bf K}\leq c$ then $\boldsymbol{\Delta}d_F(x_0,x)\geq (n-1){\bf
  ct}_c(d_F(x_0,x))$ for every $x\in M\setminus \{x_0\};$
  \item[{\rm (b)}] If $c\leq {\bf K}$ then $\boldsymbol{\Delta}d_F(x_0,x)\leq (n-1){\bf ct}_c(d_F(x_0,x))$ for every $x\in M\setminus \{x_0\}.$
\end{itemize}
\end{theorem}

\textit{2.4. Reversibility and uniformity constants.} Inspired by Rademacher
\cite{Rademacher}, we introduce the \textit{reversibility constant}
associated with $F$,
\begin{equation}  \label{reverzibilis}
r_{F}=\sup_{x\in M}r_F(x)\ \ \ \mathrm{where}\ \ \ r_F(x)=\sup_{\substack{ y
\in T_x M\setminus \{0\}}} \frac{F(x,y)}{F(x,-y)}.
\end{equation}
It is clear that $r_{F}\geq 1$ (possibly, $r_{F}=+\infty$) and $r_{F}= 1$ if
and only if $(M,F)$ is reversible. In the same way, we define the constant $%
r_{F^*}$ associated with $F^*$ and one has $r_{F^*}=r_F.$

The number
\begin{equation*}
l_{F}=\inf_{x\in M}l_F(x)\ \ \ \mathrm{where}\ \ \ l_F(x)= \inf_{y,v,w\in
T_xM\setminus \{0\}}\frac{g_{(x,v)}(y,y)}{g_{(x,w)}(y,y)},
\end{equation*}
is the \textit{uniformity constant of $F$} which measures how far $F$ and $%
F^*$ are from Riemannian structures, see Egloff \cite{Egloff}. Indeed, one
can see that $l_{F}\leq 1$, and $l_{F}= 1$ if and only if $(M,F)$ is a
Riemannian manifold, see Ohta \cite{Ohta-Math-Ann}. In the same manner, we
can define the constant $l_{F^*}$ for $F^*$, and it follows that $%
l_{F^*}=l_F.$ The definition of $l_{F}$ in turn shows that
\begin{equation}  \label{eq:2uni}
F^{*2}( x,{t\alpha+(1-t)\beta} ) \le tF^{*2}(x,\alpha) +(1-t)F^{*2}(x,\beta)
-{l_{F}}t(1-t) F^{*2}(x,\beta-\alpha)
\end{equation}
for all $x\in M$, $\alpha,\beta \in T_x^*M$ and $t\in [0,1]$.

By the above definitions, one can easily deduce that
\begin{equation}  \label{rev-unif}
l_F(x)r_F^2(x)\leq 1,\ x\in M.
\end{equation}
\noindent For the Randers metric (\ref{Randers-metrika}), a direct
computation gives that
\begin{equation}  \label{rev-unif-11}
r_F(x)=\frac{1+\|\beta\|_h(x)}{1-\|\beta\|_h(x)}\ \ \mathrm{and}\ \
l_F(x)=\left(\frac{1-\|\beta\|_h(x)}{1+\|\beta\|_h(x)}\right)^2,\
x\in M,
\end{equation}
see also Yuan and Zhao \cite{Yuan-Zhao}.

\begin{proposition}
\label{prop-rev-unifor} Let $(M,F)$ be a Finsler manifold. Then the
following statements hold:

\begin{itemize}
\item[\textrm{(a)}] If $l_F>0$ then $r_F<+\infty;$

\item[\textrm{(b)}] If $r_F<+\infty$, then forward and backward completeness
of $(M,F)$ coincide$;$

\item[\textrm{(c)}] If $(M,F)$ is of Randers type $($see {{\rm(\ref
{Randers-metrika})}$)$} with $\mathbf{S}=0$ then $l_F>0.$
\end{itemize}
\end{proposition}

\textit{Proof.} (a) follows by (\ref{rev-unif}). (b) is a simple consequence
of the Hopf-Rinow theorem, since a set in $M$ is forward bounded if and only
if it is backward bounded whenever $r_F<+\infty.$ (c) If $(M,F)$ is of
Randers type with $\mathbf{S}=0$ and $F$ has the form from (\ref%
{Randers-metrika}), Ohta \cite{Ohta-Randers} proved that $\beta$ is
a Killing form of constant $h-$length, i.e., there exists
$\beta_0\in (0,1)$ such
that $\|\beta\|_h(x)=\beta_0$ for every $x\in M$. Therefore, by (\ref%
{rev-unif-11}), one has that $l_F=\left(\frac{1-\beta_0}{1+\beta_0}%
\right)^2>0.$ \hfill $\square$

\section{Reversibility versus Sobolev spaces on non-compact Finsler
manifolds}\label{sect-3}

\textit{Proof of Theorem \ref{thm-1}.} Due to the  convexity of
$F^{*2}$, if $u,v\in W^{1,2}_0(M,F,\mathsf{m})$ then $u+v\in
W^{1,2}_0(M,F,\mathsf{m}).$ Moreover, since $r_F<\infty$, one also
has that $cu\in W^{1,2}_0(M,F,\mathsf{m})$ for every $c\in \mathbb
R$ and $u\in W^{1,2}_0(M,F,\mathsf{m})$. Consequently,
$W^{1,2}_0(M,F,\mathsf{m})$ is a vector space over $\mathbb R.$

Note that $\|\cdot\|_{F_s}$ is a norm and $\|\cdot\|_{F}$ is an
asymmetric norm. Moreover, a simple argument based on the definition
of the reversibility constant $r_F$ gives that $\|\cdot\|_{F_s}$ and
$\|\cdot\|_{F}$ are equivalent; in particular, one has
$$\left(\frac{1+r_F^{2}}{2}\right)^{-{1}/{2}}F^*(x,\alpha)\leq
F_s^*(x,\alpha)\leq
\left(\frac{1+r_F^{-2}}{2}\right)^{-{1}/{2}}F^*(x,\alpha),\ \forall
(x,\alpha)\in T^*M;$$ thus relation (\ref{norm-equivalent}) also
yields.

Let
$$L^2(M,\mathsf{m})=\left\{u:M\to \mathbb R:u\ {\rm is\ measurable},\ \|u\|_{L^2(M,\mathsf{m})}<\infty\right\},$$
where $$\|u\|_{L^2(M,\mathsf{m})}=\left(\int_M
u^2(x)\mathrm{d}\mathsf{m}(x)\right)^{1/2}.$$ It is standard that
$(L^2(M,\mathsf{m}),\|\cdot\|_{L^2(M,\mathsf{m})})$ is a Hilbert
space. Since $F^{*2}$
is a (strictly) convex function, so $F_s^{*2}$, one can prove that $%
(W_{0}^{1,2}(M,F,\mathsf{m}),\|\cdot\|_{F_s})$ is a closed subspace of the Hilbert space $L^2(M,\mathsf{m})$, which concludes the proof. \hfill $%
\square$

\begin{remark}\rm
\textrm{The statement of Theorem \ref{thm-1} remains valid for an
arbitrary open domain $\Omega\subset M$ instead of the whole
manifold $M$. }
\end{remark}

\textit{Sharpness of Theorem \ref{thm-1}.} We claim that in general
$W^{1,2}_0(M,F,\mathsf{m})$ need not has a vector space structure.
In fact,  we \textit{cannot} assert that $u\in
W^{1,2}_0(M,F,\mathsf{m})$ implies $-u\in W^{1,2}_0(M,F,\mathsf{m})$
whenever the reversibility constant $r_F$ is not finite. A similar
phenomenon is already pointed out in \cite{Kristaly-NA} for a
Funk-type metric.

For completeness, we provide another example on the
Finsler-Poincar\'e disc model. If $x_1=r\cos \theta$ and $%
x_2=r\sin \theta$ are the polar coordinates, let $$M=B^2(0,2)=\{x=(x_1,x_2)\in \mathbb{%
R}^2:r^2=|x|^2=x_1^2+x_2^2<4\}$$ and $F$ be the Randers metric given by (\ref%
{Randers-metrika}) where
\begin{equation*}
h=\frac{16}{(4-r^2)^2}(\mathrm{d}r^2+r^2\mathrm{d}\theta^2)\
\mathrm{and}\ \beta=\frac{16r}{16-r^4}\mathrm{d}r,
\end{equation*}
see \cite[Section 12.6]{BCS}.
 Consequently, if
$V=p\frac{\partial}{\partial r}+q\frac{\partial}{\partial \theta}\in
T_{(r,\theta)}M$, then we explicitly have
\begin{equation*}
F((r,\theta),V)=\frac{4}{4-r^2}\sqrt{p^2+r^2q^2}+\frac{16pr}{16-r^4}.
\end{equation*}
Note that
\begin{equation}  \label{beta-becsles}
\|\beta\|_h(x)=\frac{4r}{4+r^2},
\end{equation}
thus the volume element (see (\ref{randers-volume})) takes the form
\begin{equation}  \label{FP-volume-form}
\mathrm{d}\mathsf{m}(x)=
\mathrm{d}V_F(x)=\frac{16r(4-r^2)}{(4+r^2)^3}\mathrm{d}r\mathrm{d}\theta.
\end{equation}
The pair $(M,F)$ is a forward (but not backward) complete Randers space with
constant flag curvature $\mathbf{K}=-\frac{1}{4}$ and
\begin{equation*}
d_F(\mathbf{0},x)=\log \left(\frac{4+r^2}{(2-r)^2}\right)\ \mathrm{and}\
d_F(x,\mathbf{0})=\log \left(\frac{(2+r)^2}{4+r^2}\right),
\end{equation*}
where $\mathbf{0}=(0,0).$

Due to relations (%
\ref{rev-unif-11}) and (\ref{beta-becsles}), one has $r_F(x)=\left(\frac{2+r%
}{2-r}\right)^2,$ where $r=|x|.$ Consequently, the reversibility
constant
\begin{equation*}
r_F=\sup_{x\in M}r_F(x)=\lim_{r\to
2}\left(\frac{2+r}{2-r}\right)^2=+\infty
\end{equation*}
and $$ l_F=0.$$

By (\ref{derivalt-local}), one has
\begin{equation*}
Dd_F(\mathbf{0},x)=\frac{4(2+r)}{(2-r)(4+r^2)}\mathrm{d}r.
\end{equation*}
Therefore, by means of (\ref{Randers-dual}), a direct computation yields
that
\begin{equation}  \label{alabb-kell---}
F^*(x,Dd_F(\mathbf{0},x))=1\ \mathrm{and}\ F^*(x,-Dd_F(\mathbf{0},x))=\left(%
\frac{2+r}{2-r}\right)^2.
\end{equation}
Note that the first relation in (\ref{alabb-kell---}) follows also by (\ref%
{tavolsag-derivalt}).

Let $u:M\to \mathbb{R}$ be defined by $$u(x)=-e^{-\frac{d_F(\mathbf{0},x)}{4}%
}.$$ It is clear that $u\in W^{1,2}_{\mathrm{loc}}(M)$. Since $Du(x)=\frac{1}{%
4}e^{-\frac{d_F(\mathbf{0},x)}{4}}D d_F(\mathbf{0},x),$ by the first
relation of (\ref{alabb-kell---}) and (\ref{FP-volume-form}) one has
\begin{eqnarray*}
I_+ &:=& \int_M F^{*2}(x,Du(x))\mathrm{d}\mathsf{m}(x)=\frac{1}{16}\int_M e^{-\frac{%
d_F(\mathbf{0},x)}{2}}\mathrm{d}\mathsf{m}(x) \\
&=& 2\pi \int_0^2\frac{r(2-r)^2(2+r)}{(4+r^2)^\frac{7}{2}}\mathrm{d}r \\
&=& \frac{\pi}{30},
\end{eqnarray*}
thus $u\in W^{1,2}(M,F,\mathsf{m})$. Furthermore,
$$I:=\int_M u^2(x)\mathrm{d}\mathsf{m}(x)=\int_M e^{-\frac{%
d_F(\mathbf{0},x)}{2}}\mathrm{d}\mathsf{m}(x)=\frac{8\pi}{15}.$$
Thus, $$\|u\|_F^2=I_++I=\frac{17\pi}{30},$$ so $u\in
W_0^{1,2}(M,F,\mathsf{m})$.

 However, the second relation of (\ref{alabb-kell---}) and (\ref%
{FP-volume-form}) imply that
\begin{eqnarray*}
I_- &:=& \int_M F^{*2}(x,-Du(x))\mathrm{d}\mathsf{m}(x)=\frac{1}{16}\int_M e^{-%
\frac{d_F(\mathbf{0},x)}{2}}F^{*2}(x,-Dd_F(\mathbf{0},x))\mathrm{d}\mathsf{m}(x) \\
&=& 2\pi \int_0^2\frac{r(2+r)^5}{(4+r^2)^\frac{7}{2}(2-r)^2}\mathrm{d}r \\
&=& +\infty,
\end{eqnarray*}
Therefore,
$$\|-u\|_F^2=I_-+I=+\infty,$$
i.e., $-u\notin W^{1,2}(M,F,\mathsf{m})$ and $-u\notin
W_0^{1,2}(M,F,\mathsf{m})$. Moreover, according to
(\ref{szimmetrikus-randers}), one has that
\begin{align*}
\|u\|_{F_s}^2=\|-u\|_{F_s}^2
=&\frac{\pi}{5}+\frac{\pi}{16}\sqrt{-2+2\sqrt{2}}\ln\left(5+4\sqrt{2}+4\sqrt{4-2\sqrt{2}}+6\sqrt{-2+2\sqrt{2}}\right) \\
&+\frac{\pi}{8}\sqrt{2+2\sqrt{2}}\arctan\left(\frac{\sqrt{2+2\sqrt{2}}-\sqrt{-2+2\sqrt{2}}}{\sqrt{2}-2}\right)
\\
\approx& 0.1877.
\end{align*}
Consequently, the norm $\|\cdot\|_{F_s}$ and the asymmetric norm
$\|\cdot\|_{F}$ are not equivalent.



\section{Convexity of the singular Hardy-Finsler energy functional}

In order to deal with singular problems of type $({\mathcal{P}}%
^{\mu}_\Omega) $ we first need a Hardy inequality on (not necessarily
reversible) Finsler-Hadamard manifold with $\mathbf{S}=0$. As mentioned
before, these spaces include Finsler-Hadamard manifolds of Berwald type
(thus, both Minkowski spaces and Hadamard-Riemannian manifolds).

\begin{proposition}
\label{Hardy} Let $(M,F)$ be an $n-$di\-men\-sional $(n\geq 3)$
Finsler-Hadamard manifold with $\mathbf{S}=0$, and let $x_0\in M$ be fixed.
Then
\begin{equation}
\int_{M}F^{*2}(x,-D(|u|)(x))\mathrm{d}{\mathsf m}(x) \geq \overline \mu\int_{M}\frac{%
u^2(x)}{d_F^2(x_0,x)}\mathrm{d}{\mathsf m}(x),\ \forall u\in
C_0^\infty(M),
\end{equation}
where the constant $\overline \mu=\frac{(n-2)^2}{4}$ is optimal and never
achieved.
\end{proposition}

\textit{Proof.} By convexity, we have the following inequality
\begin{equation}  \label{fontos}
F^{*2}(x,\beta)\geq F^{*2}(x,\alpha)+2(\beta-\alpha)(J^*(x,\alpha)), \
\forall \alpha,\beta\in T_x^*M.
\end{equation}
Let $x_0\in M$ and $u\in C_0^\infty(M)$ be arbitrarily fixed and let $\gamma=%
\sqrt{\overline \mu}=\frac{n-2}{2}>0$. We consider the function $%
v(x)=d_F(x_0,x)^\gamma u(x)$. Therefore, $u(x)=d_F(x_0,x)^{-\gamma} v(x)$
and one has
\begin{equation*}
D(|u|)(x)=-\gamma
d_F(x_0,x)^{-\gamma-1}|v|Dd_F(x_0,x)+d_F(x_0,x)^{-\gamma}D(|v|)(x).
\end{equation*}
Applying the inequality (\ref{fontos}) with the choices $\beta=-D|u|$ and $%
\alpha=\gamma d_F(x_0,x)^{-\gamma-1}|v| Dd_F(x_0,x),$ respectively, one can
deduce that
\begin{eqnarray*}
F^{*2}(x,-D(|u|)(x)) &\geq& F^{*2}(x,\gamma d_F(x_0,x)^{-\gamma-1}|v(x)|
Dd_F(x_0,x)) \\
&&-2d_F(x_0,x)^{-\gamma}D(|v|)(x)(J^*(x,\gamma d_F(x_0,x)^{-\gamma-1}|v(x)|
Dd_F(x_0,x))).
\end{eqnarray*}
Due to relation (\ref{tavolsag-derivalt}), to the fact that $%
J^*(x,Dd_F(x_0,x))=\boldsymbol{\nabla}d_F(x_0,x)$ and $D(|v|)(x)\in
T_x^*M,$ we obtain
\begin{equation*}
F^{*2}(x,-D(|u|)(x))\geq \gamma^2 d_F(x_0,x)^{-2\gamma-2}|v(x)|^2-2\gamma
d_F(x_0,x)^{-2\gamma-1}|v(x)|D(|v|)(x)(\boldsymbol{\nabla}d_F(x_0,x)).
\end{equation*}
Integrating the latter inequality over $M$, it yields
\begin{equation*}
\int_{M}F^{*2}(x,-D(|u|)(x))\mathrm{d}{\mathsf m}(x) \geq
\gamma^2\int_{M}d_F(x_0,x)^{-2\gamma-2}|v(x)|^2\mathrm{d}{\mathsf
m}(x)+R_0,
\end{equation*}
where
\begin{equation*}
R_0 = -2\gamma\int_M d_F(x_0,x)^{-2\gamma-1}|v(x)|D(|v|)(x)(\boldsymbol{%
\nabla}d_F(x_0,x))\mathrm{d}{\mathsf m}(x).
\end{equation*}
Since $\mathbf{S}=0$ and $\mathbf{K}\leq 0$, Theorem
\ref{comparison-lagrange}(a) shows that
\begin{equation*}
d_F(x_0,x)\boldsymbol{\Delta}d_F(x_0,x)\geq n-1\ \mathrm{for\ a.e.}\ x\in M.
\end{equation*}
Consequently, by (\ref{Green}), (\ref{tavolsag-derivalt}) and the
latter estimate one has
\begin{eqnarray*}
R_0 &=&-\gamma \int_M D(|v|^2)(d_F(x_0,x)^{-2%
\gamma-1}\boldsymbol{\nabla}d_F(x_0,x))\mathrm{d}{\mathsf m}(x) \\
&=&\gamma \int_M
|v(x)|^2{\rm div}(d_F(x_0,x)^{-2%
\gamma-1}\boldsymbol{\nabla}d_F(x_0,x))\mathrm{d}{\mathsf m}(x) \\
&=& \gamma \int_M |v(x)|^2d_F(x_0,x)^{-2\gamma-2}\left(-2\gamma-1+d_F(x_0,x)%
\boldsymbol{\Delta}d_F(x_0,x)\right)\mathrm{d}{\mathsf m}(x)\geq 0,
\end{eqnarray*}
which completes the first part of the proof.

We now prove that $\gamma^2=\frac{(n-2)^2}{4}$ is sharp. Fix the numbers $%
R>r>0$ and a smooth cutoff function $\psi:M\to [0,1]$ with supp$(\psi)=%
\overline{B^+(x_0,R)}$ and $\psi(x)=1$ for $x\in B^+(x_0,r)$, and for every $%
\varepsilon>0,$ let $u_\varepsilon(x)=(\max\{\varepsilon,d_F(x_0,x)\})^{-%
\gamma},\ x\in M. $

One the one hand, by (\ref{tavolsag-derivalt}) we have
\begin{eqnarray*}
I_1(\varepsilon)&:=&\int_{M}F^{*2}(x,-D(\psi u_\varepsilon)(x))\mathrm{d}{\mathsf m}(x) \\
&=& \int_{B^+(x_0,r)}F^{*2}(x,-D u_\varepsilon(x))\mathrm{d}{\mathsf
m}(x) + \int_{B^+(x_0,R)\setminus B^+(x_0,r)}F^{*2}(x,-D(\psi
u_\varepsilon)(x))\mathrm{d}{\mathsf m}(x) \\
&=& \gamma^2\int_{B^+(x_0,r)\setminus
B^+(x_0,\varepsilon)}d_F(x_0,x)^{-2\gamma-2}\mathrm{d}{\mathsf
m}(x)+\tilde I_1(\varepsilon),
\end{eqnarray*}
where the quantity
\begin{equation*}
\tilde I_1(\varepsilon)=\int_{B^+(x_0,R)\setminus
B^+(x_0,r)}F^{*2}(x,-D(\psi u_\varepsilon)(x))\mathrm{d}{\mathsf
m}(x)
\end{equation*}
is finite and does not depend on $\varepsilon>0$ whenever $\varepsilon<r.$
On the other hand,
\begin{eqnarray*}
I_2(\varepsilon)&:=&\int_{M}\frac{(\psi u_\varepsilon)^2(x)}{d_F(x_0,x)^2}\mathrm{d}{\mathsf m}(x) \\
&\geq& \int_{B^+(x_0,r)\setminus B^+(x_0,\varepsilon)}d_F(x_0,x)^{-2\gamma-2}%
\mathrm{d}{\mathsf m}(x)=:\tilde I_2(\varepsilon).
\end{eqnarray*}
By (\ref{volume-comp-altalanos-0}), one has
\begin{equation*}
\mathrm{Vol}_F(B^+(x_0,\rho))\geq \omega_n \rho^n,\ \forall \rho>0.
\end{equation*}
Therefore, by applying the layer cake representation, we deduce that for $%
0<\varepsilon<r$, one has
\begin{eqnarray*}
\tilde I_2(\varepsilon)&=&\int_{B^+(x_0,r)\setminus
B^+(x_0,\varepsilon)}d_F(x_0,x)^{-2\gamma-2}\mathrm{d}{\mathsf m}(x)
= \int_{B^+(x_0,r)\setminus
B^+(x_0,\varepsilon)}d_F(x_0,x)^{-n}\mathrm{d}{\mathsf m}(x)
\\
&\geq& \int_{r^{-n}}^{\varepsilon^{-n}}\mathrm{Vol}_F(B^+(x_0,\rho^{-\frac{1}{n}%
})){\text d}\rho \\
&\geq& \omega_n\int_{r^{-n}}^{\varepsilon^{-n}}\rho^{-1}{\text d}\rho \\
&=&n\omega_n(\ln r-\ln \varepsilon).
\end{eqnarray*}
In particular, $\lim_{\varepsilon\to 0^+}\tilde I_2(\varepsilon)=+\infty.$
Thus, from the above relations it follows that
\begin{eqnarray*}
\frac{(n-2)^2}{4} &\leq & \inf_{u\in C_0^\infty(M)\setminus\{0\}}\frac{%
\displaystyle\int_{M}F^{*2}(x,-D(|u|)(x))\mathrm{d}{\mathsf m}(x)}{\displaystyle%
\int_{M}\frac{ u^2(x)}{d_F^2(x_0,x)}\mathrm{d}{\mathsf m}(x)} \\
&\leq& \lim_{\varepsilon\to
0^+}\frac{I_1(\varepsilon)}{I_2(\varepsilon)}
\leq\lim_{\varepsilon\to 0^+}\frac{\gamma^2\tilde
I_2(\varepsilon)+\tilde
I_1(\varepsilon)}{\tilde I_2(\varepsilon)} \\
&=&\gamma^2=\frac{(n-2)^2}{4}.
\end{eqnarray*}
A standard reasoning shows that this constant is never achieved. \hfill $%
\square$

\begin{remark}\rm
\textrm{Proposition \ref{Hardy} can be proved for an arbitrary open domain $%
\Omega\subset M$ instead of the whole manifold $M$ with $x_0\in \Omega.$ }
\end{remark}

In the sequel, we prove the main result of this section.

\begin{theorem}
\label{energy-convex} Let $(M,F)$ be an $n-$dimensional $(n\geq 3)$
Finsler-Hadamard manifold with $\mathbf{S}=0$ and $l_F>0.$ Let $%
\Omega\subseteq M$ be an open domain and $x_0\in \Omega$. Then the
functional  $\mathscr{K}_\mu:W_{0}^{1,2}(\Omega,F,{\mathsf m})\to
\mathbb{R}$ defined by
\begin{equation*}
\mathscr{K}_\mu(u)=\int_{\Omega}F^{*2}(x,D u(x))\mathrm{d}%
{\mathsf
m}(x)-\mu\int_{\Omega}\frac{u^2(x)}{d_F^2(x_0,x)}\mathrm{d}{\mathsf
m}(x)
\end{equation*}
is positive unless $u=0$ and strictly convex whenever $0\leq \mu <
l_{F} r_{F}^{-2} \overline \mu.$
\end{theorem}

\textit{Proof.} Let $0\leq \mu < l_{F} r_{F}^{-2} \overline \mu$ and $x_0\in
\Omega$ be fixed arbitrarily. By (\ref{rev-unif}), one has $r_F^2\leq l_F^{-1}<+\infty.$ The positivity of $\mathcal{K}%
_\mu $ follows by Proposition \ref{Hardy}. Let $0<t<1$ and $u,v\in
W_{0}^{1,2}(\Omega,F,{\mathsf m})$, $u\neq v$ be fixed. Then, by
(\ref{eq:2uni}), from the fact that
\begin{equation*}
F^*(x,D(v-u)(x))\geq r_{F}^{-1} F^*(x,-D(|v-u|)(x)),\ x\in \Omega,
\end{equation*}
and Proposition \ref{Hardy}, one has
\begin{eqnarray*}
\mathscr{K}_\mu\left(tu+(1-t)v\right) &=& \int_{\Omega}F^{*2}(x,tD
u(x)+(1-t)Dv(x))\mathrm{d}{\mathsf m}(x) -{\mu}\int_{\Omega}\frac{(tu+(1-t)v)^2}{d_F(x_0,x)^2}\mathrm{d}{\mathsf m}(x) \\
&\leq& t\int_{\Omega} F^{*2}(x,Du(x))\mathrm{d}{\mathsf m}(x)
+(1-t)\int_{\Omega}
F^{*2}(x,Dv(x))\mathrm{d}{\mathsf m}(x) \\
&& -{l_{F}}t(1-t) \int_{\Omega} F^{*2}(x,D(v-u)(x))\mathrm{d}{\mathsf m}(x) \\
&& -{\mu}\int_{\Omega}\frac{(tu+(1-t)v)^2}{d_F^2(x_0,x)}\mathrm{d}{\mathsf m}(x) \\
&=& t\mathscr{K}_\mu\left(u\right) +(1-t)\mathscr{K}_\mu\left(v\right) \\
&&- t(1-t)l_F\int_{\Omega} \left( F^{*2}(x,D(v-u)(x))-\mu l_F^{-1}\frac{%
(v-u)^2}{d_F^2(x_0,x)}\right)\mathrm{d}{\mathsf m}(x) \\
&\leq& t\mathscr{K}_\mu\left(u\right) +(1-t)\mathscr{K}_\mu\left(v\right) \\
&&- t(1-t)l_Fr_F^{-2}\int_{\Omega} \left( F^{*2}(x,-D|v-u|(x))-\mu
l_F^{-1}r_F^{2}\frac{(v-u)^2}{d_F^2(x_0,x)}\right)\mathrm{d}{\mathsf m}(x) \\
&<&t\mathscr{K}_\mu\left(u\right)
+(1-t)\mathscr{K}_\mu\left(v\right),
\end{eqnarray*}
which concludes the proof. \hfill $\square$

\section{Singular Poisson equations on Finsler-Hadamard manifolds}

\label{Section-Dirichlet}


Let $(M,F)$ be a (not necessarily reversible) complete, $n-$dimensional ($%
n\geq 3$) Finsler manifold, and $\Omega\subset M$ be an open domain,
$x_0\in \Omega$. For $\mu\in \mathbb{R}$, on
$W_{0}^{1,2}(\Omega,F,{\mathsf m})$ we define the \textit{singular
Finsler-Laplace operator}
\begin{equation*}
\mathcal{L}_F^\mu u =\boldsymbol{\Delta} (-u)-\mu\frac{u}{d_F^2(x_0,x)}.
\end{equation*}

\begin{proposition}
\label{comparison} {\rm (Comparison principle)}
Let $(M,F)$ be an $n-$dimensional $(n\geq 3)$ Finsler-Hadamard manifold with
$\mathbf{S}=0$ and $l_F>0.$ Let $\Omega\subset M$ be an open domain. If $%
\mathcal{L}_F^\mu u\leq \mathcal{L}_F^\mu v$ in $\Omega$ and {$u\leq v$ on $%
\partial \Omega,$} then $u\leq v$ a.e. in $\Omega$ whenever $\mu \in [0,
l_Fr_F^{-2}\overline \mu).$
\end{proposition}

\textit{Proof.} Assume that $\Omega_+=\{x\in \Omega:u(x)>v(x)\}$ has a
positive measure. Then, multiplying $\mathcal{L}_F^\mu u\leq \mathcal{L}%
_F^\mu v$ by $(u-v)_+$, by (\ref{Green}) one obtains
\begin{equation*}
\int_{\Omega_+} (D(-v)-D(-u))(\boldsymbol{\nabla}(-v)-\boldsymbol{\nabla}%
(-u))\mathrm{d}{\mathsf
m}(x)-\mu\int_{\Omega_+}\frac{(u-v)^2}{d_F^2(x_0,x)}\mathrm{d}{\mathsf
m}(x)\leq 0.
\end{equation*}
By (\ref{grad-deriv}) and the mean value theorem, the definition of $l_F$
yields that for every $x\in \Omega_+$,
\begin{eqnarray*}
(D(-v)-D(-u))(\boldsymbol{\nabla}(-v)-\boldsymbol{\nabla}(-u))&\geq &
l_FF^{*2}(x,D(-v)-D(-u)) \\
&=&l_FF^{*2}(x,D(u-v)) \\
&\geq & l_Fr_F^{-2} F^{*2}(x,-D(u-v)).
\end{eqnarray*}
Combining these relations with Proposition \ref{Hardy}, it follows that
\begin{equation*}
\left(l_Fr_F^{-2}-\frac{\mu}{\overline \mu}\right)\int_{%
\Omega_+}F^{*2}(x,-D(u-v)(x))\mathrm{d}{\mathsf m}(x)\leq 0,
\end{equation*}
which is a contradiction. \hfill $\square$\newline

Let $\mu\in [0,l_{F} r_{F}^{-2} \overline \mu)$ and $\kappa\in
L^\infty(\Omega)$. We consider the singular Poisson problem
\begin{equation*}
\ \left\{
\begin{array}{lll}
\mathcal{L}_F^\mu u=\kappa(x) & \mbox{in} & \Omega; \\
u= 0 & \mbox{on} & \partial \Omega,
\end{array}%
\right. \eqno{({\mathcal P}_\Omega^{\mu,\kappa})}
\end{equation*}
where $\Omega\subset M$ is an open, bounded domain. We introduce the \textit{%
singular energy functional} associated with the operator
$\mathcal{L}_F^\mu $ on $W_{0}^{1,2}(\Omega,F,{\mathsf m})$, defined
by
\begin{equation*}
\mathcal{E}_\mu(u)= (\mathcal{L}_F^\mu u)(u).
\end{equation*}
According to (\ref{Green}), we have in fact
\begin{equation*}
\mathcal{E}_\mu(u)=\int_{\Omega}F^{*2}(x,-D u(x))\mathrm{d}{\mathsf m}(x)-\mu\int_{M}%
\frac{u^2(x)}{d_F(x_0,x)^2}\mathrm{d}{\mathsf
m}(x)=\mathscr{K}_\mu(-u).
\end{equation*}

\begin{theorem}
\label{unique-theo} Let $(M,F)$ be an $n-$dimensional $(n\geq 3)$
Finsler-Hadamard manifold with $\mathbf{S}=0$ and $l_F>0.$ Let $%
\Omega\subset M$ be an open, bounded domain and a non-negative function $%
\kappa\in L^\infty(\Omega)$. Then problem $({\mathcal{P}}_\Omega^{\mu,%
\kappa})$ has a unique, non-negative weak solution for every $\mu\in
[0,l_Fr_F^{-2}\overline \mu)$.
\end{theorem}

\textit{Proof.} Let $\mu \in \lbrack 0,l_{F}r_{F}^{-2}\overline{\mu })$ be
fixed and consider the energy functional associated with problem $({\mathcal{%
P}}_{\Omega }^{\mu ,\kappa })$, i.e.,
\begin{equation*}
\mathcal{F}_{\mu }(u)=\frac{1}{2}\mathscr{K}_{\mu }(-u)-\int_{\Omega
}\kappa (x)u(x)\mathrm{d}{\mathsf m}(x),\ \ u\in
W_{0}^{1,2}(\Omega,F,{\mathsf m}).
\end{equation*}%
It is clear that $\mathcal{F}_{\mu }\in C^{1}(W_{0}^{1,2}(\Omega,F,{\mathsf m}),\mathbb{R}%
)$, and its critical points are precisely the weak solutions of problem $({%
\mathcal{P}}_{\Omega }^{\mu ,\kappa })$. Let $R>0$ and $x_{0}\in M$ be such
that $\Omega \subset B^{+}(x_{0},R).$ According to Wu and Xin \cite[Theorem
7.3]{Wu-Xin}, we have
\begin{equation*}
\lambda _{1}(\Omega )=\inf_{u\in W_{0}^{1,2}(\Omega,F,{\mathsf m})\setminus \{0\}}\frac{%
\displaystyle\int_{\Omega }F^{\ast 2}(x,Du(x))\mathrm{d}{\mathsf m}(x)}{%
\displaystyle\int_{\Omega }u^{2}(x)\mathrm{d}{\mathsf m}(x)}\geq \frac{(n-1)^{2}}{%
4R^{2}r_{F}^{2}}.
\end{equation*}%
Consequently, for every $u\in W_{0}^{1,2}(\Omega,F,{\mathsf m}),$
one has that
\begin{equation*}
\int_{\Omega }F^{\ast 2}(x,Du(x))\mathrm{d}{\mathsf m}(x)\geq
\frac{\lambda _{1}(\Omega )}{1+\lambda _{1}(\Omega )}\Vert u\Vert
_{F}^{2}.
\end{equation*}%
Since $\Vert \cdot \Vert _{F}$ and $\Vert \cdot \Vert _{F_{s}}$ are
equivalent (see (\ref{norm-equivalent})), we conclude that
$\mathcal{F}_{\mu }$ is bounded from below and coercive on the
reflexive Banach space $(W_{0}^{1,2}(\Omega,F,{\mathsf m}),\Vert
\cdot \Vert _{F_{s}})$, i.e., $\mathcal{F}_{\mu }(u)\rightarrow
+\infty $ whenever
$\Vert u\Vert _{F_{s}}\rightarrow +\infty $. Due to Theorem \ref%
{energy-convex}, $\mathcal{F}_{\mu }$ is strictly convex on $%
W_{0}^{1,2}(\Omega,F,{\mathsf m})$, thus the basic result of the
calculus of variations
implies that $\mathcal{F}_{\mu }$ has a unique (global) minimum point $%
u_{\mu }\in W_{0}^{1,2}(\Omega,F,{\mathsf m})$ of $\mathcal{F}_{\mu }$, see Zeidler \cite[%
Theorem 38.C and Proposition 38.15]{Zeidler}, which is also the
unique critical point of $\mathcal{F}_{\mu }.$ Since $\kappa \geq
0,$  Proposition \ref{comparison} implies that $u_{\mu }\geq 0.$
\hfill $\square $

\begin{remark}\rm Theorem \ref{uniq-elso} directly follows by
Theorem \ref{unique-theo}.
\end{remark}

\begin{lemma}\label{laplacef(d)} Let $f\in C^{2}(0,\infty)$ be a non-increasing
function. Then $$ \mathcal{L}_F^\mu (f(d_{F}(x_{0},x)))=-f^{\prime
\prime
}(d_{F}(x_{0},x))-f^{\prime }(d_{F}(x_{0},x))\cdot \mathbf{\Delta }%
d_{F}(x_{0},x)-\mu \frac{f(d_{F}(x_{0},x))}{d_{F}^{2}(x_{0},x)},\
x\in M\setminus \{x_0\}. $$
\end{lemma}

\textit{Proof.} Since $f'\leq 0,$ the claim follows from basic
properties of the Legendre transform. Namely, one has
\begin{eqnarray*}
  \mathbf{\Delta }%
(-f(d_{F}(x_{0},x))) &=& \mathrm{div}(\boldsymbol{\nabla }(-f(d_{F}(x_{0},x))))=\mathrm{div}\left( J^{\ast }(x,D(-f(d_{F}(x_{0},x)))\right) \\
   &=& \mathrm{div}(J^{\ast }(x,-f^{\prime
}(d_{F}(x_{0},x))Dd_{F}(x_{0},x))) \\
   &=& \mathrm{div}(-f^{\prime }(d_{F}(x_{0},x))\boldsymbol{\nabla }%
d_{F}(x_{0},x))\\
&=& -f^{\prime \prime }(d_{F}(x_{0},x))-f^{\prime
}(d_{F}(x_{0},x))\cdot \mathbf{\Delta }d_{F}(x_{0},x),
\end{eqnarray*}
which concludes the proof. \hfill $\square $\\

For every  $\mu \in \lbrack 0, \overline{\mu })$, $c\leq 0$ and
$\rho>0$, we recall the
 ordinary differential equation
\begin{equation*}
\left\{
\begin{array}{ll}
f^{\prime \prime }(r)+(n-1) f^{\prime }(r) {\bf ct}_{c}(r)+\mu
\frac{f(r)}{r^{2}}+1=0, \ r\in (0,\rho ], \\
f(\rho )=0,\ \displaystyle\int_0^\rho f'(r)^2r^{n-1}{\rm d}r<\infty.
\end{array}%
\right. \eqno{({\mathcal Q}^{\mu}_{c,\rho})}
\end{equation*}

\begin{proposition} \label{prop-ODE}
 $({\mathcal Q}^{\mu}_{c,\rho})$ has a unique, non-negative,
non-increasing solution belonging to $C^\infty(0,\rho)$.
\end{proposition}

{\it Proof.} We fix  $\mu \in \lbrack 0, \overline{\mu })$, $c\leq
0$ and $\rho>0$. Let us consider the Riemannian space form $(M,g_c)$
with  constant sectional curvature $c\leq 0$, i.e., $(M,g_c)$ is
isometric to the Euclidean space when $c=0$, or $(M,g_c)$ is
isometric to the hyperbolic space with sectional curvature $c<0.$
Let $x_0\in M$ be fixed. Since $(M,g_c)$ verifies the assumptions of
Theorem \ref{unique-theo}, problem
\begin{equation*}
\ \left\{
\begin{array}{lll}
-{\Delta}_{g_c} u-\mu\frac{u}{d_{g_c}^2(x_0,x)}=1 & \mbox{in} & B_{g_c}(x_0,\rho); \\
u= 0 & \mbox{on} & \partial B_{g_c}(x_0,\rho),
\end{array}%
\right. \eqno{({\mathcal R}_{c,\rho}^{\mu})}
\end{equation*}
has a unique, non-negative solution $u_0$ which is nothing but the
unique global minimum point of the energy functional
$\mathcal{F}_{\mu }:W_{0}^{1,2}(B_{g_c}(x_0,\rho),g_c,{\mathsf
m})\to \mathbb R$ defined by
\begin{equation*} \mathcal{F}_{\mu }(u)=\frac{1}{2}\int_{B_{g_c}(x_0,\rho)}|D u(x)|_{g_c}^2\mathrm{d}%
{\mathsf
m}(x)-\frac{\mu}{2}\int_{B_{g_c}(x_0,\rho)}\frac{u^2(x)}{d_{g_c}^2(x_0,x)}\mathrm{d}{\mathsf
m}(x)-\int_{B_{g_c}(x_0,\rho) }u(x)\mathrm{d}{\mathsf m}(x).
\end{equation*}
In this particular case, ${\rm d}{\mathsf m}$ denotes the canonical
Riemannian volume form on $(M,g_c)$.

Let $u_0^\star:B_{g_c}(x_0,\rho)\to [0,\infty)$ be the
non-increasing symmetric rearrangement of $u_0$ in the space form
$(M,g_c)$, see Baernstein \cite{Baernstein}. Note that
P\'{o}lya-Szeg\H{o} and Hardy-Littlewood inequalities imply that
$$\int_{B_{g_c}(x_0,\rho)}|D u_0(x)|_{g_c}^2\mathrm{d}%
{\mathsf m}(x)\geq\int_{B_{g_c}(x_0,\rho)}|D u_0^\star(x)|_{g_c}^2\mathrm{d}%
{\mathsf m}(x),$$ and
$$\int_{B_{g_c}(x_0,\rho)}\frac{u_0^2(x)}{d_{g_c}^2(x_0,x)}\mathrm{d}{\mathsf
m}(x)\leq
\int_{B_{g_c}(x_0,\rho)}\frac{{u_0^\star}^2(x)}{d_{g_c}^2(x_0,x)}\mathrm{d}{\mathsf
m}(x),$$ respectively. Moreover, by the Cavalieri principle, we also
have that
$$\int_{B_{g_c}(x_0,\rho) }u_0(x)\mathrm{d}{\mathsf m}(x)=\int_{B_{g_c}(x_0,\rho) }u_0^\star(x)\mathrm{d}{\mathsf m}(x).$$
Therefore, we  obtain that $\mathcal{F}_{\mu }(u_0)\geq
\mathcal{F}_{\mu }(u_0^{\star }).$ Consequently, by the uniqueness
of the global minimizer of $\mathcal{F}_{\mu }$ we have
$u_0=u_0^{\star }$; thus, its form is $u_0(x)=f(t)$ where
$t=d_{g_c}(x_0,x)$ and $f:(0,\rho]\to \mathbb R$ is a non-negative
and non-increasing function.  Clearly, $f(\rho)=0$ since $u_0(x)=0$
whenever $d_{g_c}(x_0,x)=\rho$.
 Moreover, since $u_0=u_0^{\star }\in
W_{0}^{1,2}(B_{g_c}(x_0,\rho),g_c,{\mathsf m})$, a suitable change
of variables gives that $\displaystyle\int_0^\rho f'(r)^2r^{n-1}{\rm
d}r<\infty$. By Lemma \ref{laplacef(d)} and Theorem
\ref{comparison-lagrange} it follows that the first part of
$({\mathcal R}_{c,\rho}^{\mu})$ can be transformed into the first
part of $({\mathcal Q}_{c,\rho}^{\mu})$; in particular, problem
$({\mathcal Q}_{c,\rho}^{\mu})$ has a non-negative, non-increasing
solution. Standard regularity theory implies that $f\in
C^\infty(0,\rho),$ see Evans \cite[p. 334]{Evans}. Finally, if we
assume that $({\mathcal Q}_{c,\rho}^{\mu})$ has two distinct
non-negative, non-increasing solutions $f_1$ and $f_2$, then both
functions $u_i(x)=f_i(d_{g_c}(x_0,x))$ ($i\in \{1,2\}$) verify
$({\mathcal R}_{c,\rho}^{\mu})$, which are distinct global minima of
the energy functional $\mathcal{F}_{\mu }$, a contradiction.
 \hfill $\square$\\


\textit{Proof of Theorem \ref{estimate_with_constant_flag}.} Let $u$ be the unique solution of problem $({%
\mathcal{P}}_{\Omega}^{\mu })$.   We claim that
\begin{equation*}
\left\{
\begin{array}{lll}
\mathcal{L}_F^\mu(\sigma _{\mu ,\rho_1 ,c_{1}}(d_{F}(x_{0},x)))\leq
1=\mathcal{L}_F^\mu(u)& {\rm
in}& B^{+}(x_{0},\rho_1 ); \\
\sigma _{\mu ,\rho_1 ,c_{1}}(d_{F}(x_{0},x))=0\leq u(x) & {\rm on} &
\partial B^{+}(x_{0},\rho_1 ),
\end{array}%
\right.
\end{equation*}
where $\rho_1=\sup\{\rho>0:B^+(x_0,\rho)\subset \Omega\}$. On one
hand, since $c_1\leq {\bf K}$, due to Theorem
\ref{comparison-lagrange} (b) and to the fact that $\sigma _{\mu
,\rho_1 ,c_{1}}$ is non-increasing, by equation $({\mathcal
Q}^{\mu}_{c_1,\rho_1})$ one has for $x\in B^{+}(x_{0},\rho_1
)\setminus \{x_0\}$,
\begin{eqnarray*}
  1 &=& -\sigma _{\mu
,\rho_1 ,c_{1}}^{\prime \prime }(d_{F}(x_{0},x))-(n-1) \sigma _{\mu
,\rho_1 ,c_{1}}^{\prime }(d_{F}(x_{0},x)) {\bf
ct}_{c_1}(d_{F}(x_{0},x))-\mu  \frac{\sigma _{\mu
,\rho_1 ,c_{1}}(d_{F}(x_{0},x))}{d_{F}^2(x_{0},x)} \\
   &\geq& -\sigma _{\mu
,\rho_1 ,c_{1}}^{\prime \prime }(d_{F}(x_{0},x))- \sigma _{\mu
,\rho_1
,c_{1}}^{\prime }(d_{F}(x_{0},x))\mathbf{\Delta }%
d_{F}(x_{0},x)-\mu  \frac{\sigma _{\mu
,\rho_1 ,c_{1}}(d_{F}(x_{0},x))}{d_{F}^2(x_{0},x)} \\
   &=&\mathcal{L}_F^\mu(\sigma _{\mu ,\rho_1 ,c_{1}}(d_{F}(x_{0},x))).
\end{eqnarray*}
 On the other
hand, since $u$ is non-negative in $\Omega$, it follows that
$0=\sigma _{\mu ,\rho_1 ,c_{1}}(d_{F}(x_{0},x))\leq u(x)$ on
$\partial B^{+}(x_{0},\rho_1 ).$ It remains to apply the comparision
principle (Proposition \ref{comparison}), obtaining $$\sigma _{\mu
,\rho_1 ,c_{1}}(d_{F}(x_{0},x))\leq u(x)\ {\rm for\ a.e.}\ x\in
B^{+}(x_{0},\rho_1 ).$$

Similarly, by using Theorem \ref{comparison-lagrange} (a) and $ {\bf
K}\leq c_2$, one can prove that
\begin{equation*}
\left\{
\begin{array}{lll}
1=\mathcal{L}_F^\mu(u)\leq \mathcal{L}_F^\mu(\sigma _{\mu ,\rho_2
,c_{2}}(d_{F}(x_{0},x)))& {\rm
in}& \Omega; \\
u(x)=0\leq \sigma _{\mu ,\rho_2 ,c_{2}}(d_{F}(x_{0},x)) & {\rm on} &
\partial \Omega,
\end{array}%
\right.
\end{equation*}
where $\rho_2=\inf\{\rho>0:\Omega\subset B^+(x_0,\rho) \}.$ In
particular, by Proposition \ref{comparison} again we have that
$$ u(x)\leq\sigma _{\mu
,\rho_2 ,c_{2}}(d_{F}(x_{0},x))\ {\rm for\ a.e.}\ x\in \Omega.$$

If $\mathbf{K}=c\leq 0$ and $\Omega=B^+(x_0,\rho)$ for some
$\rho>0,$ then $\rho_1=\rho_2=\rho$, and from above it follows that
$u(x)=\sigma _{\mu ,\rho
,c}(d_{F}(x_{0},x))$ is the unique weak solution of problem $({%
\mathcal{P}}_{B^{+}(x_{0},\rho )}^{\mu })$ which is also a pointwise
solution in $B^+(x_0,\rho)\setminus \{x_0\}.$ \hfill $\square$\\

A simple consequence of Theorem \ref{estimate_with_constant_flag} is
the following

\begin{corollary}
\label{Minkowski-proposition} Let $(M,F)=(\mathbb{R}^n,\|\cdot\|)$
be a Minkowski space and let $\mu\in
[0,l_Fr_F^{-2}\overline \mu),$ $x_0\in \mathbb{R}^n$ and $\rho>0$ be fixed. Then $%
u=\sigma_{\mu,\rho,0}(\|\cdot-x_0\|)\in
C^\infty(B^+(x_0,\rho)\setminus\{x_0\})$ is the unique pointwise
solution to problem $({\mathcal{P}}^{\mu}_{B^+(x_0,\rho)})$ in
$B^+(x_0,\rho)\setminus\{x_0\}.$



\end{corollary}

\textit{Proof.}  $(M,F)=(\mathbb{R}^{n},\Vert \cdot \Vert )$ being a
Minkowski space, it is a Finsler-Hadamard manifold with $\mathbf{S}=0$, ${\bf K}=0$ and $%
l_{F}>0.$ It remains to apply Theorem
\ref{estimate_with_constant_flag}.
\hfill$%
\square $

\begin{remark}\rm
\textrm{(i) In addition to the conclusions of Corollary
\ref{Minkowski-proposition}, one can see that }

\begin{itemize}
\item[(a)] \textrm{$\sigma_{\mu,\rho,0} \in C^1(B^+(x_0,\rho))$ if and only if $%
\mu=0$, and }

\item[(b)] \textrm{$\sigma_{\mu,\rho,0} \in C^2(B^+(x_0,\rho))$ if and only if $%
\mu=0$ and $F=\|\cdot\|$ is Euclidean. }
\end{itemize}

\textrm{(ii) When $(M,F)=(\mathbb{R}^n,\|\cdot\|)$ is a reversible
Minkowski space and $\mu=0$, Corollary \ref{Minkowski-proposition}
reduces to Theorem 2.1 from Ferone and Kawohl \cite{FeKa}. }
\end{remark}

In connection with Corollary \ref{Minkowski-proposition} we
establish an estimate for the solution of the singular Poisson
equation on \textit{backward} geodesic balls on Minkowski spaces. To
do this, we assume that  $\sigma_{\mu,r_F^{-1}\rho,0}$ is extended
beyond $r_F^{-1}\rho$ formally by the same function, its explicit form being given after the problem ${({\mathcal Q}^{\mu}_{c,\rho})}.$ Although problem $({%
\mathcal{P}}^{\mu}_{B^-(x_0,\rho)})$ cannot be solved explicitly in
general, the following sharp estimates can be given for its unique
solution by means of the reversibility constant $r_F$.

\begin{proposition}
\label{Minkowski-proposition-2} Let $(M,F)=(\mathbb{R}^n,\|\cdot\|)$
be a Minkowski space and let $\mu\in [0,l_Fr_F^{-2}\overline \mu),$
$x_0\in \mathbb{R}^n$ and $\rho>0$ be fixed. If $\tilde
u_{\mu,\rho}$ denotes the unique weak solution to problem
$({\mathcal{P}}^{\mu}_{B^-(x_0,\rho)})$, then
\begin{equation*}
(\sigma_{\mu,r_F^{-1}\rho,0}(\|x-x_0\|))_+\leq\tilde
u_{\mu,\rho}(x)\leq \sigma_{\mu,r_F\rho,0}(\|x-x_0\|)\ { for\ a.e.}\
x\in B^-(x_0,\rho).
\end{equation*}
  Moreover, the above two bounds coincide if and
only if $(M,F)$ is reversible.
\end{proposition}

\textit{Proof.} The proof immediately follows by the comparison
principle Proposition \ref{comparison}, showing that
\begin{equation*}
\ \left\{
\begin{array}{lll}
\mathcal{L}_F^{\mu}(w_{\mu,\rho}^{-})= 1 = \mathcal{L}_F^{\mu}(w_{\mu,%
\rho}^{+}) & \mbox{in} & B^-(x_0,\rho); \\
w_{\mu,\rho}^{-}\leq 0 \leq w_{\mu,\rho}^{+} & \mbox{on} & \partial
B^-(x_0,\rho),
\end{array}%
\right.
\end{equation*}
where $w_{\mu,\rho}^{-}(x)=\sigma_{\mu,r_F^{-1}\rho,0}(\|x-x_0\|)$ and $%
w_{\mu,\rho}^{+}(x)=\sigma_{\mu,r_F\rho,0}(\|x-x_0\|)$, respectively. \hfill $%
\square$\\

\textit{Proof of Theorem \ref{rigidity}.}  Let $x_0\in M$ be fixed
and we assume that for some $\mu\in [0,l_{F} r_{F}^{-2} \overline
\mu)$, the function $u(x)=\sigma _{\mu ,\rho ,c}(d_{F}(x_{0},x))$ is
the unique pointwise solution of $({\mathcal{P}}_{B^{+}(x_{0},\rho
)}^{\mu })$ on $B^+(x_0,\rho)\setminus \{x_0\}$ for some $\rho
>0.$  By Lemma \ref{laplacef(d)}
and from the fact that $\sigma _{\mu ,\rho ,c}$ is a solution of
$({\mathcal Q}^{\mu}_{c,\rho})$, it follows that
$$\boldsymbol{\Delta}d_F(x_0,x)= (n-1){\bf ct}_c(d_F(x_0,x))\ {\rm
in}\ B^{+}(x_{0},\rho )\setminus \{x_0\}$$ pointwisely. The latter
relation and a simple calculation shows that
\begin{equation}\label{delta-egyenloseg}\mathbf{\Delta }%
w_c(d_{F}(x_{0},x))=1\ {\rm in}\ B^{+}(x_{0},\rho )\setminus
\{x_0\},$$ where
$$w_c(r)=\int_{0}^{r}{\bf s}_c(s)^{-n+1}\int_{0}^{s}{\bf s}_c(t)%
^{n-1}{\rm d}t{\rm d}s.\end{equation} Let $0<\tau<\rho$ be fixed
arbitrarily. The unit outward normal vector to the forward geodesic
sphere $S^+(x_0,\tau)=\partial B^+(x_0,\tau)=\{x\in
M:d_F(x_0,x)=\tau\}$ at $x\in S^+(x_0,\tau)$ is given by $\mathbf{n}=%
\boldsymbol{\nabla} d_F(x_0,x)$. Let us denote by ${\text
d}\varsigma_F(x)$ the canonical volume form on $S^+(x_0,\tau)$
induced by $\mathrm{d}{\mathsf m}(x)={\text d}V_F(x)$.
By Stokes' formula (see \cite{Shen-Annalen}, \cite[Lemma 3.2]{Wu-Xin}) and $%
g_{(x,\mathbf{n})}(\mathbf{n},\mathbf{n})=F(x,\mathbf{n})^2=F(x, \boldsymbol{%
\nabla} d_F(x_0,x))^2=1$ (see (\ref{tavolsag-derivalt})), on account
of relation  (\ref{delta-egyenloseg})
 we have
\begin{eqnarray*}
\mathrm{Vol}_{F}(B^{+}(x_{0},\tau )) &=&\int_{B^{+}(x_{0},\tau )}\boldsymbol{%
\Delta }(w_c(d_F(x_0,x)))\mathrm{d}{\mathsf
m}(x) \\
&=&\int_{B^{+}(x_{0},\tau )}\mathrm{div}(\boldsymbol{\nabla
}(w_c(d_F(x_0,x))))\mathrm{d}{\mathsf
m}(x) \\
&{=}&\int_{S^{+}(x_{0},\tau )}g_{(x,\mathbf{n})}(\mathbf{n},w_c^{\prime }(d_F(x_0,x))%
\boldsymbol{\nabla }d_F(x_0,x)){\text{d}}\varsigma _{F}(x) \\
&=&w_c^{\prime }(\tau )\cdot \mathrm{Area}_{F}({S^{+}(x_{0},\tau
)}).
\end{eqnarray*}%
Therefore,
\begin{equation*}
\frac{\mathrm{Area}_{F}({S^{+}(x_{0},\tau )})}{\mathrm{Vol}%
_{F}(B^{+}(x_{0},\tau ))}=\frac{1}{w_c^{\prime }(\tau )}=\frac{{\bf s}_c%
(\tau )^{n-1}}{\displaystyle\int_{0}^{\tau }{\bf s}_c%
(t)^{n-1}{\rm d}t},
\end{equation*}%
or equivalently,
\begin{equation*}
\frac{\frac{\rm d}{\rm d\tau }\mathrm{Vol}_{F}(B^{+}(x_{0},\tau ))}{\mathrm{Vol}%
_{F}(B^{+}(x_{0},\tau ))}=\frac{\frac{\rm d}{\rm d\tau }\displaystyle\int_{0}^{\tau }{\bf s}_c%
(t)^{n-1}{\rm d}t}{\displaystyle\int_{0}^{\tau }{\bf s}_c%
(t)^{n-1}{\rm d}t}.
\end{equation*}%
Integrating the latter expression on the interval $[s,\rho ],$
$0<s<\rho,$ and exploiting (\ref{volume-comp-nullaban}), it follows
that
\begin{equation}\label{egyenloseg-utolssso}
\frac{\mathrm{Vol}_{F}(B^{+}(x_{0},\rho
))}{V_{c,n}(\rho)}=\lim_{s\to
0^+}\frac{\mathrm{Vol}_{F}(B^{+}(x_{0},s ))}{V_{c,n}(s)}=1.
\end{equation}
%
%
 According to Theorem \ref{comparison-volume}, it yields
$${\bf K}(\cdot;\dot \gamma_{x_0,y}(t))=c$$ for every $t\in
[0,\rho)$ and $y\in T_{x_0}M$ with $F(x_0,y)=1$, where
$\gamma_{x_0,y}$ is the constant speed geodesic with
$\gamma_{x_0,y}(0)=x_0$ and $\dot
\gamma_{x_0,y}(0)=y.$ This concludes the proof. \hfill$%
\square $\newline

\textit{Proof of Theorem \ref{rigid_Minkowski}}. The implications
"(a)$\Rightarrow $(b)" and  "(c)$\Rightarrow $(a)" are trivial, see
Corollary \ref{Minkowski-proposition}; it remains to prove
"(b)$\Rightarrow $(c)". By the proof of Theorem \ref{rigidity} we
know the validity of relation (\ref{egyenloseg-utolssso}) for every
$\rho>0$. Let $x\in M$ and $\rho>0$ be arbitrarily fixed. We have
that
\begin{eqnarray*}
1 &\leq&\frac{{\rm Vol}_F(B^+(x,\rho))}{V_{0,n}(\rho)}\ \ \ \ \ \ \
\ \ \ \
 \ \ \ \ \ \ \ \ \ \ \ \ \ \ \ \ \ \
 \ \ \ \ \ \ \ \ \ \ \ \ \ \ \ \ \ \ \ \ \ \ \ {\rm (see\
(\ref{volume-comp-altalanos-0}))} \\
&\leq&\limsup_{r\to \infty}\frac{{\rm Vol}_F(B^+(x,r))}{V_{0,n}(r)}
\ \ \ \ \ \ \ \ \ \ \ \ \ \ \ \ \
 \ \ \ \ \ \ \ \ \ \ \ \ \ \ \ \ \ \
 \ \ \ \ \ \ \ {\rm (monotonicity\ from\ Theorem\ \ref{comparison-volume})}\\&\leq& \limsup_{r\to
\infty}\frac{{\rm Vol}_F(B^+(x_0,r+d_F(x_0,x)))}{V_{0,n}(r)}\ \ \ \
\ \ \ \ \ \  \ \ \ \ \ \ \ \ \ \ \ \ \ \ \ (B^+(x,r)\subset
B^+(x_0,r+d_F(x_0,x))) \\
   &=& \limsup_{r\to
\infty}\left(\frac{{\rm Vol}_F(B^+(x_0,r+d_F(x_0,x)))}{V_{0,n}(r+d_F(x_0,x))}\cdot\frac{V_{0,n}(r+d_F(x_0,x))}{V_{0,n}(r)}\right) \\
   &=&1, \ \ \ \ \ \ \ \ \ \ \
 \ \ \ \ \ \ \ \ \ \ \ \ \ \ \ \ \ \
 \ \ \ \ \ \ \ \ \ \ \ \ \ \ \ \ \ \ \ \ \ \ \ \ \ \ \ \ \ \ \ \  \ \ \ \ \  \ \ \ \ \ {\rm
(see\ (\ref{egyenloseg-utolssso}))}
\end{eqnarray*}
because one has
\begin{equation}\label{kkk}
\lim_{r\to
\infty}\frac{V_{0,n}(r+d_F(x_0,x))}{V_{0,n}(r)}=\lim_{r\to
\infty}\frac{(r+d_F(x_0,x))^n}{r^n}=1.
\end{equation}
Consequently,
\begin{equation}\label{vol-id-neg-1}
    {{\rm Vol}_F(B^+(x,\rho))}= V_{0,n}(\rho)=\omega_n\rho^n\ {\rm for\ all}\ x\in M\ {\rm
and}\ \rho>0.
\end{equation}
On account of Theorem \ref{comparison-volume} and relation
(\ref{vol-id-neg-1}), we conclude that
 ${\bf K}=0$.

Note that every Berwald space with ${\bf K}=0$ is necessarily a
locally Minkowski space, see Bao, Chern and Shen
\cite[Section~10.5]{BCS}. Therefore, the global volume identity
  (\ref{vol-id-neg-1}) actually implies that $(M,F)$ is isometric to a Minkowski
space.
 \hfill $\square$

\vspace{1cm}

\noindent \textbf{Acknowledgment.} Research supported by a grant of
the Romanian National Authority for Scientific Research,
CNCS-UEFISCDI, project no. PN-II-ID-PCE-2011-3-0241. Cs. Farkas is
also supported by the State of Hungary, co-financed by the European
Social Fund in the framework of T\'AMOP 4.2.4.A/2-11-1-2012-0001
'National Excellence Program' and by Collegium Talentum. The
research of A. Krist\'aly is also supported by J\'anos Bolyai
Research Scholarship.


\end{document}